\documentclass[11pt]{article}
\usepackage{amsmath,epsfig}
\newcommand{\comment}[1]{}

\newcommand{\fata}{\mathbf{a}}
\newcommand{\fatb}{\mathbf{b}}
\newcommand{\fatc}{\mathbf{c}}
\newcommand{\fatd}{\mathbf{d}}
\newcommand{\reals}{\mbox{$\mathbb R$}}

\newcommand{\nats}{\mbox{$\mathbb N$}}
\newcommand{\ints}{\mbox{$\mathbb Z$}}

\newcommand{\floor}[1]{\lfloor #1 \rfloor}
\newcommand{\ceil}[1]{\lceil #1 \rceil}

\newcommand{\conv}{{\mathop{\mathrm{conv}}\nolimits}}
\newcommand{\Pri}{{\mathop{\mathrm{Pri}}\nolimits}}
\newcommand{\Pyr}{{\mathop{\mathrm{Pyr}}\nolimits}}
\newcommand{\supp}{{\mathop{\mathrm{supp}}\nolimits}} 

\usepackage{amssymb}

\def\squarebox#1{\hbox to #1{\hfill\vbox to #1{\vfill}}}
\def\qed{\hspace*{\fill}
        \vbox{\hrule\hbox{\vrule\squarebox{.667em}\vrule}\hrule}\smallskip}
\newenvironment{proof}{\begin{trivlist}
  \item[\hspace{\labelsep}{\em\noindent Proof.~}]
  }{\qed\end{trivlist}}

\newtheorem{lemma}{Lemma}[section]
\newtheorem{theorem}[lemma]{Theorem}
\newtheorem{corollary}[lemma]{Corollary}
\newtheorem{proposition}[lemma]{Proposition}
\newtheorem{claim}[lemma]{Claim}
\newtheorem{observation}[lemma]{Observation}
\newtheorem{definition}[lemma]{Definition}

\def\squareforqed{\hbox{\rlap{$\sqcap$}$\sqcup$}}
\def\qed{\ifmmode\squareforqed\else{\unskip\nobreak\hfil
\penalty50\hskip1em\null\nobreak\hfil\squareforqed
\parfillskip=0pt\finalhyphendemerits=0\endgraf}\fi}

\newlength{\tablength}
\newlength{\spacelength}
\newcommand{\tabstar}{\hspace*{\tablength}}
\newcommand{\spacestar}{\hspace*{\spacelength}}
\def\obeytabs{\catcode`\^^I=\active}
{\obeytabs\global\let^^I=\tabstar}
{\obeyspaces\global\let =\spacestar}
\newenvironment{display}{\begingroup\obeylines\obeyspaces\obeytabs}{\endgroup}
\newenvironment{prog}{\begin{display}\parskip0pt\sf}{\end{display}}

\title{The flag polynomial of the Minkowski sum of simplices}

\author{
{\sl Geir Agnarsson} 
\thanks{Department of Mathematical Sciences,
George Mason University,
MS 3F2,
4400 University Drive,
Fairfax, VA -- 22030, USA,
{\tt geir@math.gmu.edu}}
}

\date{}

\begin{document}

\maketitle

\begin{abstract}
For a polytope we define the {\em flag polynomial},
a polynomial in commuting variables related to 
the well-known flag vector and describe 
how to express the the flag polynomial 
of the Minkowski sum of $k$ standard simplices
in a direct and canonical way in terms of the
{\em $k$-th master polytope} $P(k)$ where $k\in\nats$.
The flag polynomial facilitates many direct computations.
To demonstrate this we provide two examples; we first derive a  
formula for the $f$-polynomial and the maximum number
of $d$-dimensional faces of the Minkowski sum of two simplices.
We then compute the maximum discrepancy between the number
of $(0,d)$-chains of faces of a Minkowski sum of two simplices
and the number of such chains of faces of a simple polytope
of the same dimension and on the same number of vertices.

\vspace{3 mm}



\vspace{2 mm}

\end{abstract}

\section{Introduction and motivation}
\label{sec:intro}

For a polytope $P$ of dimension $d$, the {\em face vector} of $P$ 
is given by $f(P) = (f_0(P),\ldots, f_{d-1}(P))$ where
$f_i(P)$ is the number of $i$-dimensional faces of $P$. It is an ongoing
investigation to describe the face vectors for $d$-dimensional convex 
polytopes. The case $d\leq 3$ are the only cases where this description
is complete. A nice theorem by Steinitz~\cite{Steinitz} states that an integer
vector $(f_0,f_1,f_2)$ is a face vector for a three dimensional polytope $P$ if and only if
(i) we have the Euler relation $f_0 - f_1 + f_2 = 2$, (ii) $f_2 - 4 \leq 2(f_0 - 4)$
and (iii) $f_0 - 4 \leq 2(f_2 - 4)$. Moreover, the first inequality is tight if and only
if $P$ is simplicial, the second one is tight if and only if $P$ is simple. The only 
linear equation satisfied by the entries of the face vector is the Euler-Poincar\'{e} 
relation 
$\sum_{i=0}^{d-1}(-1)^if_i = 1 + (-1)^{d-1}$, that holds for all polytopes of dimension $d$.
For higher dimensions $d\geq 4$ no such nice description of face vectors
has been found yet. In the hope of obtaining a description of a polytope 
in terms of parameters/entries of a vector, there is the {\em flag vector}
$(f_S(P))_{S\subseteq \{0,1,\ldots,d-1\}}$, containing $2^d$ entries where
for each $S = \{s_1,s_2,\ldots,s_{\ell}\}\subseteq \{0,\ldots,d-1\}$ the
entry $f_S(P)$ denotes the number of strict chains (or {\em flags}) 
$A_1\subset A_2\subset\cdots\subset A_{\ell}$ of faces of $P$ where $\dim(A_i) = s_i$. 
The entries $f_S(P)$ do satisfy the generalized Dehn-Sommerville equations
which, together with the existence of matching polytopes, implies that the linear 
span of flag vectors of $d$-dimensional polytopes is precisely $F_d$, 
the $d$-th Fibonacci number, as shown in the seminal paper~\cite{Bayer-Billera}. 
One way to encode the information contained in the flag vector
of a $d$-polytope $P$ is to define the $\fata\fatb$-index, a linear combination 
of non-commuting
monomials in $\fata$ and $\fatb$ of degree $d$ in the following way: for 
$S\subseteq \{0,1,\ldots,d-1\}$
let $w_S = t_0t_1\cdots t_{d-1}$, where $t_i = \fata$ if $i\not\in S$ and 
$t_i = \fatb$ if $i\in S$. 
The {\em $\fata\fatb$-index} $\Psi(P)$ for $P$ is then defined as 
$\Psi(P) = \sum_S h_S(P) w_S$,
where the sum is over all $S\subseteq \{0,1,\ldots,d-1\}$ and 
$h_S(P) = \sum_{T\subseteq S}(-1)^{|S\setminus T|}f_T(P)$.
An elegant consequence of the generalized Dehn-Sommerville equations for the 
entries $f_S(P)$ is that each $\fata\fatb$-index can be written as a 
non-commuting polynomial in $\fatc = \fata+\fatb$ and $\fatd = \fata\fatb+\fatb\fata$.
The resulting polynomial called the {\em $\fatc\fatd$-index}, also 
denoted $\Psi(P)$ and first
appeared in~\cite{Bayer-Klapper}.
Within the free non-commuting $\reals$-algebra $\reals\langle \fata,\fatb\rangle$ in 
two variables $\fata$ and $\fatb$, the degree of $\fatc$ is set to be 1 and the 
degree of $\fatd$ is set to be two. 
In this way the number of all homogeneous monomials of degree
$d$ in $\fatc$ and $\fatd$ is precisely $F_d$. Hence, the $\fatc\fatd$-index is the 
most compact
way to encode all the entries of the flag vector of a polytope $P$ w.r.t.~linear
relations of the coefficients of the monomial in $\fatc$ and $\fatd$. A lot of current
research is therefore directed toward understanding the $\fatc\fatd$-index for polytopes.
Although the $\fatc\fatd$-index is an elegant and compact way to encode the flag vector
for a polytope $P$, it can be at times difficult to handle in many explicit computations since
algebraic manipulation in 
$\reals\langle \fatc,\fatd\rangle\subseteq \reals\langle \fata,\fatb\rangle$ 
can become particularly cumbersome.

The purpose of this article is to do the opposite of the above: namely, instead
of attempting to compactify the encoding of the flag vector further, then we shall
instead relax the constrictions and allow some redundancies among coefficients in 
order to facilitate direct computations. This will in certain cases allow us to 
write down explicit (albeit involved)
formulae from which the entries of the flag vector can be obtained. 
We will define the {\em flag polynomial} for a polytope,
a polynomial in commuting variables that will behave reasonably well when 
constructing polytopes
by known methods i.e.~prisms, pyramids, products, free joins and Minkowski sums.

The article is organized as follows: In Section~\ref{sec:basic} we present our basic
definitions and write some direct consequences. In Section~\ref{sec:flag-poly} we state
(but do not prove) Theorem~\ref{thm:flag-main}, the main result of the section; a formula
for the flag polynomial of a Minkowski sum of $k$ standard simplices. We then parse
through an explicit application of Theorem~\ref{thm:flag-main} and note some
consequences. Section~\ref{sec:join} is devoted to the derivation of 
Theorem~\ref{thm:P-acts-flag}; a formula for the flag polynomial of a polytope
$P(Q_1,\ldots,Q_m)$, a polytope obtained when the geometric polytope $P$
in ${\reals}^m$ acts on the polytopes $Q_1,\ldots,Q_m$. This action 
generalizes the notion of the free join of two polytopes. Finally, we note
that Theorem~\ref{thm:flag-main} is a corollary of Theorem~\ref{thm:P-acts-flag}.
In the last Section~\ref{sec:f-poly-2} we present two examples on the application
of the flag polynomial. In the first example we compute the maximum number of
faces of an arbitrary Minkowski sum of two simplices. In the latter example we
compute the maximum discrepancy between the number of $(0,d)$-chains of faces
of a Minkowski sum of two simplices and the number of $(0,d)$-chains of faces
of a simple polytope of the same dimension and on the same number of vertices
as the Minkowski sum.

This article can be 
viewed as a self contained continuation of~\cite{Minkowski}. For related
references and history of the significance of Minkowski sums of simplices
we refer to the introduction of~\cite{Minkowski}. Also, for self-containment
we will present Definition~\ref{def:master} and some terminology from~\cite{Minkowski}.

For $r\in\nats$ we let $[r] = \{1,2,\ldots,r\}$. The {\em standard simplex}
$\Delta_{r-1} = \Delta_{[r]}$ of dimension $r-1$ is given by
$
\Delta_{r-1} = \Delta_{[r]} = 
\{ \tilde{x} = (x_1,\ldots,x_r)\in {\reals}^r : x_i\geq 0 \mbox{ for all $i$ },
x_1+\cdots+x_r = 1\}.
$
Each subset $F\subseteq [r]$ yields a {\em face} $\Delta_F$ of $\Delta_{[r]}$ given by
$
\Delta_F = \{ \tilde{x} \in \Delta_{[r]} : x_i = 0 \mbox{ for } i\not\in F\}.
$
Clearly $\Delta_F$ is itself a simplex embedded in ${\reals}^r$. If ${\cal{F}}$ is
a family of subsets of $[r]$, then we can form the {\em Minkowski sum} of simplices
\[
P_{\cal{F}} = \sum_{F\in{\cal{F}}}\Delta_F =
\left\{\sum_{F\in{\cal{F}}}\tilde{x}_F : 
\tilde{x}_F \in \Delta_F \mbox{ for each } F\in{\cal{F}}\right\}.
\]
Note that a face of any polytope $P\subseteq {\reals}^r$ (in particular of $\Delta_{[r]}$) 
is given by the set of points that maximize a linear
functional $c^T\tilde{x}$ for $\tilde{x}\in P$.

Recall the $f$-polynomial of a polytope: let $P$ be a $d$-dimensional polytope.
For each $i\in \{0,1,\ldots,d\}$ let $\mathbf{f}_i(P)$ denote the set 
of all $i$-dimensional faces of $P$ and so $f_i(P) = |\mathbf{f}_i(P)|$ is the number of 
$i$-dimensional faces of $P$. Then $\mathbf{f}(P) = \bigcup_{i=0}^{d}\mathbf{f}_i(P)$ is 
the set of all the faces of $P$ except the empty face. 
The {\em $f$-polynomial} of $P$ is then given by $f_P(x) = \sum_{i=0}^{d}f_i(P)x^i$.
\begin{observation}
\label{obs:f-poly-standard}
The $f$-polynomial for the standard $(r-1)$-dimensional simplex 
$\Delta_{r-1} = \Delta_{[r]}$ is given by
\[
f_{\Delta_{r-1}}(x) = \frac{(x+1)^r - 1}{x}.
\]
\end{observation}

\section{Definitions and some basic properties}
\label{sec:basic}

We now define our main object for this paper, a generalization
of the $f$-polynomial of a polytope, and derive some basic results.

\noindent{\sc Conventions:} (i) For an $\ell$-tuple $\tilde{x} = (x_1,\ldots,x_{\ell})$
of variables and an $\ell$-tuple of numbers $\tilde{a} = (a_1,\ldots,a_{\ell})$, 
let ${\tilde{x}}^{\tilde{a}} = x_1^{a_1}x_2^{a_2}\cdots x_{\ell}^{a_{\ell}}$.
(ii) For $\tilde{a} = (a_1,\ldots,a_{\ell})$ let
$\partial(\tilde{a}) = (a_1, a_2-a_1,a_3-a_2,\ldots,a_{\ell}-a_{\ell - 1})$.
\begin{definition}
\label{def:flag-poly}
Let $P$ be a polytope with $\dim(P) = d$ and $\ell \in \nats$. 
For an $\ell$-tuple of variables $\tilde{x} = (x_1,\ldots,x_{\ell})$  
the {\em $\ell$-flag polynomial} is defined by
\[
{\tilde{f}}^{\ell}_{P}(x_1,\ldots,x_{\ell}) = 
\sum_{\tilde{s}}f_{\tilde{s}}(P){\tilde{x}}^{\partial({\tilde{s}})},
\]
where the sum is taken over all chains $\tilde{s} = (s_1,\ldots,s_{\ell})$
with $0\leq s_1\leq s_2\leq \cdots \leq s_{\ell}\leq d$
and $f_{\tilde{s}}(P)$ denotes the number of chains of faces
$A_1\subseteq A_2\subseteq \cdots \subseteq A_{\ell}$ 
of $P$ with $\dim(A_i) = s_i$ for each $i\in\{1,\ldots, \ell\}$.
\end{definition}
Note that each $\ell$-tuple $\tilde{a}$ of non-negative integers
yields an $\ell$-tuple $\tilde{s}$ of increasing chain of integers
$s_1\leq s_2\leq \cdots\leq s_{\ell}$ where $s_i = a_1 +\cdots + a_i$
for each $i$ and so $\partial{\tilde{s}} = \tilde{a}$

The reasons why it might be worth working with the flag
polynomial of a polytope is (i) its expression is simple
for many basic polytopes like the simplex and hypercube and 
(ii) it behaves well under various natural constructions, the first of which is
the Cartesian product.
\begin{proposition}
\label{prp:prod-good}
For polytopes $P$ and $Q$, a positive integer $\ell\in \nats$ and 
$\tilde{x} = (x_1,\ldots,x_{\ell})$ we have 
${\tilde{f}}^{\ell}_{P\times Q}(\tilde{x})= 
{\tilde{f}}^{\ell}_{P}(\tilde{x}){\tilde{f}}^{\ell}_{Q}(\tilde{x})$.
\end{proposition}
\begin{proof}
Any pair of chains
$A_1\subseteq A_2\subseteq \cdots \subseteq A_{\ell}\subseteq P$ and 
$B_1\subseteq B_2\subseteq \cdots \subseteq B_{\ell}\subseteq Q$
yields a chain
\[
A_1\times B_1\subseteq A_2\times B_2
\subseteq \cdots \subseteq A_{\ell}\times B_{\ell}\subseteq P\times Q
\]
of faces in the product $P\times Q$. Vice versa, any such chain
of faces in $P\times Q$ of length $\ell$ yields a pair of chains
of faces of length $\ell$ in $P$ and $Q$ respectively. Since the 
mapping $\tilde{s}\mapsto \partial{\tilde{s}}$ is linear the proposition
follows.
\end{proof}

Consider now the standard $(r-1)$-dimensional simplex $\Delta_{r-1}$.
Let $0\leq s_1\leq s_2\leq \cdots \leq s_{\ell}\leq r-1$ be an 
$\ell$-chain of nonnegative integers. Since any choice of 
$h+1$ vertices from the $r$ vertices of $\Delta_{r-1}$ is 
an $h$-face of $\Delta_{r-1}$, the number of chains
$A_1\subseteq \cdots \subseteq A_{\ell}$ of faces of $\Delta_{r-1}$
with $\dim(A_i) = s_i$ for each $i$ is, by the product rule, given by 
\[
\binom{r}{s_{\ell}+1}\binom{s_{\ell}+1}{s_{\ell-1}+1}\cdots\binom{s_2+1}{s_1+1} = 
\binom{r}{r-s_{\ell}-1 \ \ s_{\ell}-s_{\ell -1} \ \ \cdots \ \ s_2-s_1 \ \ s_1+1}.
\]
Since $\partial(\tilde{s}) = (s_1,s_2-s_1,\ldots,s_{\ell}-s_{\ell-1})$ 
we have the following.
\begin{observation}
\label{obs:flag-poly-simp}
For $\ell,r \in\nats$ the $\ell$-flag polynomial of the standard 
simplex $\Delta_{r-1}$ is given by 
\[
{\tilde{f}}^{\ell}_{\Delta_{r-1}}(x_1,\ldots,x_{\ell}) = 
\frac{(x_1 + \cdots + x_{\ell} + 1)^r - (x_2 + \cdots + x_{\ell} + 1)^r}{x_1}.
\]
\end{observation}
By Proposition~\ref{prp:prod-good} and Observation~\ref{obs:flag-poly-simp}
we have obtain the following.
\begin{corollary}
\label{cor:flag-poly-prism}
For $\ell\in\nats$ and a polytope $P$ the $\ell$-flag polynomial
of the prism $\Pri(P) = P\times\Delta_1$ of $P$ is given by
\[
{\tilde{f}}^{\ell}_{\Pri(P)}(x_1,\ldots,x_{\ell}) = 
{\tilde{f}}^{\ell}_{P}(x_1,\ldots,x_{\ell})(2 + x_1 + 2(x_2 + \cdots + x_{\ell})).
\]
\end{corollary}
Since the $d$-dimensional hypercube $Q_d$ is given by 
$Q_d = {\Delta_1}^d = \Delta_1\times\cdots\times\Delta_1$ ($d$-times), 
then we also obtain the following.
\begin{corollary}
\label{cor:flag-poly-hcube}
For $\ell,d\in\nats$ the $\ell$-flag polynomial of the hypercube $Q_d$
is given by 
\[
{\tilde{f}}^{\ell}_{Q_d}(x_1,\ldots,x_{\ell}) = (2 + x_1 + 2(x_2 + \cdots + x_{\ell}))^d.
\]
\end{corollary}
From this we see that the number of $\ell$-chains of faces of $Q_d$ is $(2\ell+1)^d$.
In particular, for $\ell = 1$ we see that the total number of faces of $Q_d$ is
$3^d$, the conjectured maximum number of faces of a centrally symmetric
$d$-polytope~\cite{Gil-Kalai}.

Should one want to see the chains from the monomials
more explicitly we have the following directly by Definition~\ref{def:flag-poly}
of the flag polynomial.
\begin{observation}
\label{obs:part->seq}
Let $P$ be a polytope and $\ell\in\nats$. If $\tilde{x} = (x_1,\ldots, x_{\ell})$
and $\tilde{y} = (y_1,\ldots,y_{\ell})$ are two tuples of variables such that
$x_i = y_iy_{i+1}\cdots y_{\ell}$ for each $i=1,\ldots,\ell$, then
the coefficient of $y_1^{s_1}y_2^{s_2}\cdots y_{\ell}^{s_{\ell}}$
in ${\tilde{\phi}}_P^{\ell}(\tilde{y}) := {\tilde{f}}^{\ell}_{P}(\tilde{x})$
is the number of chains of faces
$A_1\subseteq A_2\subseteq \cdots \subseteq A_{\ell}$ of $P$
where $\dim(A_i) = s_i$ for each $i$.
\end{observation}

\noindent{\sc Remarks:}  
(i) If for some reason one prefers to consider a polynomial 
directly corresponding to the flag vector (i.e.~terms of the polynomial
corresponding to strict chains of faces), then that can be
obtained by the Inclusion/Exclusion principle: for each fixed $\ell$,
the polynomial corresponding to the number of {\em strict}
chains $A_1\subset \cdots \subset A_{\ell}$ of faces of a polynomial $P$
with $0\leq \dim(A_1) < \cdots < \dim(A_{\ell})\leq \dim(P)$
is given by ${\tilde{g}}^{\ell}_{P}(\tilde{x})$ where 
\[
{\tilde{g}}^{\ell}_{P}(x_1,\ldots,x_{\ell}) = 
\sum_{h=0}^{\ell - 1}
\left(    
(-1)^h\sum_{2\leq i_1<\cdots<i_h\leq\ell}
{\tilde{f}}^{\ell-h}_{P}(\widehat{x_{i_1}},\ldots,\widehat{x_{i_h}})
\right),
\]
where ${\tilde{f}}^{\ell-h}_{P}(\widehat{x_{i_1}},\ldots,\widehat{x_{i_h}})$ is 
the $(\ell-h)$-flag polynomial on the $(\ell-h)$-tuple obtained from 
$(x_1,\ldots,x_{\ell})$ by removing the variables $x_{i_j}$ for $j=1,\ldots, h$.

(ii) The number of terms in ${\tilde{g}}^{\ell}(\tilde{x})$ is $\binom{\dim(P) + 1}{\ell}$
whereas the number of terms in ${\tilde{f}}^{\ell}_{P}(\tilde{x})$ equals
the number of solutions $(a_1,\ldots,a_{\ell})$ to $a_1+\cdots+a_{\ell}\leq d = \dim(P)$
in nonnegative integers $a_i\geq 0$, namely $\binom{d+\ell}{\ell}$. Comparing
the number of entries of the flag vector to the combined number of terms of all 
flag polynomials
${\tilde{f}}^{\ell}_{P}(\tilde{x})$ for $\ell = 0,1,\ldots, d$, the combined number of terms
equals
\[
\sum_{\ell=0}^{d}\binom{d+\ell}{\ell} = \binom{2d+1}{d+1}\sim \frac{2^{2d+1}}{\sqrt{e\pi d}},
\]
not surprisingly an exponential expression which is considerably less than 
$(2^d)^2 = 4^d$,  the square of the number of entries in the flag vector,
but also more than $F_d$, the $d$-th Fibonacci number,
the exact number of linearly independent terms~\cite{Bayer-Billera}.

\vspace{3 mm}

In the next section we state Theorem~\ref{thm:flag-main}, our
main result, then parse through an elaborate example, before proving
the main result in Section~\ref{sec:join}.

\section{The flag polynomial of Minkowski sum of simplices}
\label{sec:flag-poly}

If ${\cal{F}}$ consists of $k$ subsets of $[r]$ and is presented (merely for
convenience) as an ordered tuple 
${\cal{F}} = (F_1,F_2,\ldots,F_k)$,
then for each $i \in [r]$ we let $N_{\cal{F}}(i) = \{j \in [k]:i \in F_j\}$. 
The following definition is from~\cite{Minkowski}:
\begin{definition}
\label{def:master} 
For $k\in\nats$ let 
${\cal{H}}(k) = (H_1,\ldots, H_k)$ be the family of $k$ ordered 
subsets of $[2^k-1]$ so that for 
$i \in \{1,2,\ldots,2^k-1\}$ then $N_{{\cal{H}}(k)}(i)$ is the $i$-th 
nonempty subset of $[k]$ in the lexicographic order. 
Then the Minkowski sum of these $k$ simplices $P(k):=P_{{\cal{H}}(k)}$ is called the 
$k^{th}$ {\em master polytope}.
\end{definition}

\noindent{\sc Example:}
For $k = 2$ we have $[2^k-1] = [3] = \{1,2,3\}$ and 
the lexicographic ordering of the nonempty subsets of $[2] = \{1,2\}$ is
given by $\{1,2\} > \{1\} > \{2\}$. By definition we must have
${\cal{H}}(k) = {\cal{H}}(2) = (H_1,H_2)$ where 
$H_1 =\{1,2\}$ and $H_2 = \{1,3\}$. In this case 
$P(2) = \Delta_{\{1,2\}} + \Delta_{\{1,3\}}$ is the 2-dimensional 
rhombus lying in the plane $x + y + z = 2$ in ${\reals}^3$ that is 
the convex hull of the points 
$(2,0,0)$, 
$(1,1,0)$,
$(1,0,1)$ 
and $(0,1,1)$.
In a similar fashion, for $k = 3$ we have $[2^k-1] = [7]$ and the lexicographic
ordering of the nonempty subsets of $[3]$ is given by
$\{1,2,3\} > \{1,2\} > \{1,3\} > \{1\} > \{2,3\} > \{2\} > \{3\}$. By definition
${\cal{H}}(k) = {\cal{H}}(3) = (H_1,H_2,H_3)$ where 
$H_1 = \{1,2,3,4\}$, $H_2 = \{1,2,5,6\}$ and $H_3 = \{1,3,5,7\}$.
In this case $P(3) = \Delta_{H_1} + \Delta_{H_2} + \Delta_{H_3}$ 
is a 6-dimensional polytope in ${\reals}^7$ with 41 vertices
as verified by the computer program POLYMAKE~\cite{Gawrilow}
(see~\cite[Observations 2.13 and 2.16]{Minkowski}). 

\vspace{3 mm}

Note that for a general $k\in\nats$ then $P(k) = \sum_{i = 1}^k\Delta_{H_i}$
where $|H_i| = 2^{k-1}$ for each $i$.
The master polytopes $P(k)$ for $k\geq 2$ are fundamental in understanding
the overall structure of a general Minkowski sum of standard simplices, as we will see in 
Theorem~\ref{thm:flag-main}, the main result in this section, 
where we express the flag polynomial of a general Minkowski sum of simplices. 

\noindent{\sc Example:} Since $P(2)$, the 2-dimensional rhombus, has the combinatorial structure 
of $\Delta_1\times \Delta_1$ we get by Proposition~\ref{prp:prod-good} and 
Observation~\ref{obs:flag-poly-simp} the following.
\begin{eqnarray*}
{\tilde{f}}^3_{P(2)}(x_1,x_2,x_3) & = & {\tilde{f}}^3_{\Delta_1\times \Delta_1}(x_1,x_2,x_3) \\
                                  & = & [{\tilde{f}}^3_{\Delta_1}(x_1,x_2,x_3)]^2 \\
                                  & = & \left(\frac{(x_1+x_2+x_3+1)^2 - 
(x_2+x_3+1)^2}{x_1}\right)^2 \\
                                  & = & (x_1+2x_2+2x_3+2)^2,
\end{eqnarray*}
and hence by Observation~\ref{obs:part->seq}
\begin{eqnarray*}
{\tilde{\phi}}^3_{P(2)}(y_1,y_2,y_3) & = & (y_1y_2y_3 + 2y_2y_3 + 2y_3 + 2)^2 \\
                                     & = & y_1^2y_2^2y_3^2 + 4y_1y_2^2y_3^2 + 
4y_1y_2y_3^2 + 4y_2^2y_3^2 + \\ 
                                     &   & 4y_1y_2y_3 + 8y_2y_3^2 + 8y_2y_3 + 
4y_3^2 + 8y_3 + 4,
\end{eqnarray*}
something which is evident by inspection.

\vspace{3 mm}

From now on, we will assume $k\in\nats$ fixed and we will always present our collection of
subsets as an ordered $k$-tuple ${\cal{F}} = (F_1,F_2,\ldots,F_k)$ of $[r]$,
and the Minkowski sum as
\begin{equation}
\label{eqn:Minkssum}
P_{\cal{F}} = \sum_{F\in{\cal{F}}}\Delta_F = \Delta_{F_1} + \cdots + \Delta_{F_k}.
\end{equation}
We will assume that ${\cal{F}}$ covers all of $[r]$, that is $\bigcup_{F\in{\cal{F}}}F = [r]$,
since otherwise we can simply restrict to this union.
For a general Minkowski sum $P_{\cal{F}}$ of simplices as in (\ref{eqn:Minkssum}),
then for each $i\in [2^k-1]$ we let $R_i = \{j\in [r]: N_{\cal{F}}(j) = N_{{\cal{H}}(k)}(i)\}$
and $r_i = |R_i|$, so $r_i$ is the number of elements in the set 
$\left(\bigcap_{\ell\in N_{{\cal{H}}(k)}(i)}F_{\ell}\right)\cap
\left(\bigcap_{\ell\not\in N_{{\cal{H}}(k)}(i)}\overline{F_{\ell}}\right)$.
Since ${\cal{F}}$ covers $[r]$ we have $r = r_1 + \cdots + r_{2^k-1}$.

Note that for any polytope $P$ we have that ${\tilde{f}}^{\ell}_{P}(x_1,\ldots,x_{\ell-1},0)$
only has terms from ${\tilde{f}}^{\ell}_{P}(x_1,\ldots,x_{\ell})$ where the coefficient 
for $x_{\ell}$ is zero. Hence, ${\tilde{f}}^{\ell}_{P}(x_1,\ldots,x_{\ell-1},0)$ is
the generating function for the number of chains 
$Q_1\subseteq \cdots\subseteq Q_{\ell-1} = Q_{\ell}$
of faces of $P$, which is equivalent to the generating function for the number of
chains  $Q_1\subseteq \cdots\subseteq Q_{\ell-1}$ of faces of $P$.
\begin{observation}
\label{obs:x_ell=0}
For a polytope $P$ we have
${\tilde{f}}^{\ell}_{P}(x_1,\ldots,x_{\ell-1},0) 
= {\tilde{f}}^{\ell-1}_{P}(x_1,\ldots,x_{\ell-1})$.
\end{observation}
\noindent{\sc Example:} Considering $P(2)$ again, we have from the above example that
\[
{\tilde{f}}^2_{P(2)}(x_1,x_2) = {\tilde{f}}^3_{P(2)}(x_1,x_2,0) = (x_1+2x_2+2)^2.
\]
Note that the total number of 2-chains is the sum of the coefficients of
${\tilde{f}}^2_{P(2)}(x_1,x_2)$, 25 altogether.

\vspace{3 mm}

\noindent{\sc Conventions:} (i) For a set of point $S\subseteq {\reals}^r$ let 
$\supp(S) = \{i \in [r] : x_i\neq 0\mbox{ for some } \tilde{x} = (x_1,\ldots,x_r) \in S\}$ 
be the {\em support} of $S$. (ii) For a face $A$ of $P$ let $\dim(A)$ 
denote the dimension of the face $A$. (iii) For a polytope $P$ 
we let ${\mathbf{C}}_{\ell}(P)$  denote the set of all 
$\ell$-chains $\tilde{A} = (A_1,\ldots,A_{\ell})$ of faces
of $P$ with $A_1\subseteq A_2\subseteq \cdots \subseteq A_{\ell}$.
(iv) For each such $\tilde{A}$ let $\dim(\tilde{A}) = (\dim(A_1),\ldots,\dim(A_{\ell}))$
denote the chain of their corresponding dimensions. 

\vspace{3 mm}

With the above conventions in mind we have 
\[
\partial(\dim(\tilde{A})) 
= (\dim(A_1),\dim(A_2)-\dim(A_1),\ldots, \dim(A_{\ell})-\dim(A_{\ell-1})).
\]
We can now state the following main result in this section, 
whose proof is in Section~\ref{sec:join}.
\begin{theorem}
\label{thm:flag-main}
For $\ell\in\nats$ and $\tilde{x} = (x_1,\ldots,x_{\ell})$ the flag $\ell$-polynomial for a 
general Minkowski sum of simplices 
$P_{\cal{F}} = \sum_{F\in{\cal{F}}}\Delta_F = \Delta_{F_1} + \cdots + \Delta_{F_k}$
is given by 
\[
{\tilde{f}}^{\ell}_{P_{\cal{F}}}(x_1,\ldots,x_{\ell}) = 
\sum_{\tilde{A}\in {\mathbf{C}}_{\ell}(P(k))}
\left(
{\tilde{x}}^{\partial(\dim(\tilde{A}))}
\prod_{i=1}^{\ell}
\left(
\prod_{{\alpha}_i \in\supp(A_i)\setminus \supp(A_{i-1})}
{\tilde{f}}^{\ell-i+1}_{\Delta_{r_{{\alpha}_i}-1}}(x_i,\ldots,x_{\ell})
\right)
\right).
\]
In the formula above we interpret $\supp(A_0) = \emptyset$, since there is no $A_0$.
\end{theorem}

For each $\tilde{A}\in {\mathbf{C}}_{\ell}(P(k))$ let 
\[
p_{\tilde{A}}(x_1,\ldots,x_{\ell}) = 
{\tilde{x}}^{\partial(\dim(\tilde{A}))}
\prod_{i=1}^{\ell}
\left(
\prod_{{\alpha}_i \in\supp(A_i)\setminus \supp(A_{i-1})}
{\tilde{f}}^{\ell-i+1}_{\Delta_{r_{{\alpha}_i}-1}}(x_i,\ldots,x_{\ell})
\right).
\]
For a concrete example we will parse through the very special case of 
Theorem~\ref{thm:flag-main}
when $k=2$ and $\ell = 2$. In this case $P_{\cal{F}} = \Delta_{F_1} + \Delta_{F_2}$ 
and $r_1 = |F_1\cap F_2|$, $r_2 = |F_1\setminus F_2|$ and $r_3 = |F_2\setminus F_1|$.
We will label the vertices of $P(2)$ as follows:
$A = (2,0,0)$, $B = (1,1,0)$, $C = (1,0,1)$ and $D = (0,1,1)$. A 2-chain of
faces of $P(2)$ will be given by $\tilde{A} = (A_1,A_2)$ where 
each face is presented as a string of the vertices whose convex combination
forms the face. As noted in the last example there are 25 2-chains of $P(2)$
and they can be grouped together by what $\dim(\tilde{A}) = (\dim(A_1),\dim(A_2))$
is as follows:
\begin{equation}
\label{eqn:2-chain}
\begin{array}{lll}
\dim(\tilde{A}) = (0,0) & : & (A,A), (B,B), (C,C), (D,D), \\
\dim(\tilde{A}) = (1,1) & : & (AB,AB), (AC,AC), (BD,BD), (CD,CD), \\
\dim(\tilde{A}) = (2,2) & : & (ABCD,ABCD), \\
\dim(\tilde{A}) = (0,1) & : & (A,AB), (A,AC), (B,AB), (B,BD), (C,AC), 
(C,CD), (D,BD), (D,CD), \\
\dim(\tilde{A}) = (1,2) & : & (AB, ABCD), (AC,ABCD), (BD,ABCD), (CD,ABCD),\\
\dim(\tilde{A}) = (0,2) & : & (A,ABCD), (B,ABCD), (C,ABCD), (D,ABCD).    
\end{array}
\end{equation}

For $X\in\{A,B,C,D\}$ and $\tilde{A} = (A_1,A_2) = (X,X)$ we have
\[
\partial(\dim(\tilde{A})) = (\dim(X),\dim(X)-\dim(X)) = (0,0),
\]
$\supp(A_1) = \supp(X)$, $\supp(A_2)\setminus\supp(A_1) = \emptyset$, and hence
\[
p_{\tilde{A}}(x_1,x_2) = p_{(X,X)}(x_1,x_2) =
\prod_{{\alpha}_1 \in\supp(X)}{\tilde{f}}^2_{\Delta_{r_{{\alpha}_1}-1}}(x_1,x_2).
\]
From this we obtain
\begin{eqnarray*}
p_{(A,A)}(x_1,x_2) & = & {\tilde{f}}^2_{\Delta_{r_1-1}}(x_1,x_2), \\
p_{(B,B)}(x_1,x_2) & = & 
{\tilde{f}}^2_{\Delta_{r_1-1}}(x_1,x_2){\tilde{f}}^2_{\Delta_{r_2-1}}(x_1,x_2), \\
p_{(C,C)}(x_1,x_2) & = & 
{\tilde{f}}^2_{\Delta_{r_1-1}}(x_1,x_2){\tilde{f}}^2_{\Delta_{r_3-1}}(x_1,x_2), \\
p_{(D,D)}(x_1,x_2) & = & 
{\tilde{f}}^2_{\Delta_{r_2-1}}(x_1,x_2){\tilde{f}}^2_{\Delta_{r_3-1}}(x_1,x_2).
\end{eqnarray*}

For $X,Y\in\{A,B,C,D\}$, and $\tilde{A} = (A_1,A_2) = (XY,XY)$
we have 
\[
\partial(\dim(\tilde{A})) = (\dim(XY),\dim(XY)-\dim(XY)) = (1,0),
\]
$\supp(A_1) = \supp(XY)$, $\supp(A_2)\setminus\supp(A_1) = \emptyset$,
and hence
\[
p_{\tilde{A}}(x_1,x_2) = p_{(XY,XY)}(x_1,x_2) = 
x_1\prod_{{\alpha}_1 \in\supp(XY)}{\tilde{f}}^2_{\Delta_{r_{{\alpha}_1}-1}}(x_1,x_2).
\]
From this we obtain
\begin{eqnarray*}
p_{(AB,AB)}(x_1,x_2) & = & 
x_1{\tilde{f}}^2_{\Delta_{r_1-1}}(x_1,x_2){\tilde{f}}^2_{\Delta_{r_2-1}}(x_1,x_2), \\
p_{(AC,AC)}(x_1,x_2) & = & 
x_1{\tilde{f}}^2_{\Delta_{r_1-1}}(x_1,x_2){\tilde{f}}^2_{\Delta_{r_3-1}}(x_1,x_2), \\
p_{(BD,BD)}(x_1,x_2) & = & x_1{\tilde{f}}^2_{\Delta_{r_1-1}}(x_1,x_2)
{\tilde{f}}^2_{\Delta_{r_2-1}}(x_1,x_2){\tilde{f}}^2_{\Delta_{r_3-1}}(x_1,x_2), \\
p_{(CD,CD)}(x_1,x_2) & = & x_1{\tilde{f}}^2_{\Delta_{r_1-1}}(x_1,x_2)
{\tilde{f}}^2_{\Delta_{r_2-1}}(x_1,x_2){\tilde{f}}^2_{\Delta_{r_3-1}}(x_1,x_2).
\end{eqnarray*}

For $\tilde{A} = (A_1,A_2) = (ABCD,ABCD)$
we have 
\[
\partial(\dim(\tilde{A})) = (\dim(ABCD),\dim(ABCD)-\dim(ABCD)) = (2,0),
\]
$\supp(A_1) = \supp(ABCD) = \{1,2,3\}$, $\supp(A_2)\setminus\supp(A_1) = \emptyset$,
and hence
\[
p_{\tilde{A}}(x_1,x_2) = p_{(ABCD,ABCD)}(x_1,x_2) = 
x_1^2\prod_{{\alpha}_1 \in\supp(ABCD)}{\tilde{f}}^2_{\Delta_{r_{{\alpha}_1}-1}}(x_1,x_2).
\]
From this we obtain
\[
p_{(ABCD,ABCD)}(x_1,x_2) = x_1^2{\tilde{f}}^2_{\Delta_{r_1-1}}(x_1,x_2)
{\tilde{f}}^2_{\Delta_{r_2-1}}(x_1,x_2){\tilde{f}}^2_{\Delta_{r_3-1}}(x_1,x_2).
\]

For $X,Y\in\{A,B,C,D\}$, and $\tilde{A} = (A_1,A_2) = (X,XY)$ we have
\[
\partial(\dim(\tilde{A})) = (\dim(X),\dim(XY)-\dim(X)) = (0,1),
\]
$\supp(A_1) = \supp(X)$, $\supp(A_2)\setminus\supp(A_1) = \supp(XY)\setminus\supp(X)$, 
and hence
\[
p_{\tilde{A}}(x_1,x_2) = p_{(X,XY)}(x_1,x_2) = 
x_2\prod_{{\alpha}_1 \in\supp(X)}{\tilde{f}}^2_{\Delta_{r_{{\alpha}_1}-1}}(x_1,x_2)
\prod_{{\alpha}_2 \in\supp(XY)\setminus\supp(X)}{\tilde{f}}^1_{\Delta_{r_{{\alpha}_2}-1}}(x_2).
\]
From this we obtain
\begin{eqnarray*}
p_{(A,AB)}(x_1,x_2) & = & 
x_2{\tilde{f}}^2_{\Delta_{r_1-1}}(x_1,x_2){\tilde{f}}^1_{\Delta_{r_2-1}}(x_2), \\
p_{(A,AC)}(x_1,x_2) & = & 
x_2{\tilde{f}}^2_{\Delta_{r_1-1}}(x_1,x_2){\tilde{f}}^1_{\Delta_{r_3-1}}(x_2), \\
p_{(B,AB)}(x_1,x_2) & = & 
x_2{\tilde{f}}^2_{\Delta_{r_1-1}}(x_1,x_2){\tilde{f}}^2_{\Delta_{r_2-1}}(x_1,x_2), \\
p_{(B,BD)}(x_1,x_2) & = & x_2{\tilde{f}}^2_{\Delta_{r_1-1}}(x_1,x_2)
{\tilde{f}}^2_{\Delta_{r_2-1}}(x_1,x_2){\tilde{f}}^1_{\Delta_{r_3-1}}(x_2), \\
p_{(C,AC)}(x_1,x_2) & = & 
x_2{\tilde{f}}^2_{\Delta_{r_1-1}}(x_1,x_2){\tilde{f}}^2_{\Delta_{r_3-1}}(x_1,x_2), \\
p_{(C,CD)}(x_1,x_2) & = & x_2{\tilde{f}}^2_{\Delta_{r_1-1}}(x_1,x_2)
{\tilde{f}}^2_{\Delta_{r_3-1}}(x_1,x_2){\tilde{f}}^1_{\Delta_{r_2-1}}(x_2), \\
p_{(D,BD)}(x_1,x_2) & = & x_2{\tilde{f}}^2_{\Delta_{r_2-1}}(x_1,x_2)
{\tilde{f}}^2_{\Delta_{r_3-1}}(x_1,x_2){\tilde{f}}^1_{\Delta_{r_1-1}}(x_2), \\
p_{(D,CD)}(x_1,x_2) & = & x_2{\tilde{f}}^2_{\Delta_{r_2-1}}(x_1,x_2)
{\tilde{f}}^2_{\Delta_{r_3-1}}(x_1,x_2){\tilde{f}}^1_{\Delta_{r_1-1}}(x_2).
\end{eqnarray*}

For $XY\in\{AB,AC,BD,CD\}$, and $\tilde{A} = (A_1,A_2) = (XY,ABCD)$ we have
\[
\partial(\dim(\tilde{A})) = (\dim(XY),\dim(ABCD)-\dim(XY)) = (1,1),
\]
$\supp(A_1) = \supp(XY)$, 
$\supp(A_2)\setminus\supp(A_1) = 
\supp(ABCD)\setminus\supp(XY) =\{1,2,3\}\setminus\supp(XY)$
and hence
\[
p_{\tilde{A}}(x_1,x_2) = p_{(XY,ABCD)}(x_1,x_2) = 
x_1x_2\prod_{{\alpha}_1 \in\supp(XY)}{\tilde{f}}^2_{\Delta_{r_{{\alpha}_1}-1}}(x_1,x_2)
\prod_{{\alpha}_2 \in\{1,2,3\}\setminus\supp(XY)}{\tilde{f}}^1_{\Delta_{r_{{\alpha}_2}-1}}(x_2).
\]
From this we obtain
\begin{eqnarray*}
p_{(AB,ABCD)}(x_1,x_2) & = & x_1x_2{\tilde{f}}^2_{\Delta_{r_1-1}}(x_1,x_2)
{\tilde{f}}^2_{\Delta_{r_2-1}}(x_1,x_2){\tilde{f}}^1_{\Delta_{r_3-1}}(x_2), \\
p_{(AC,ABCD)}(x_1,x_2) & = & x_1x_2{\tilde{f}}^2_{\Delta_{r_1-1}}(x_1,x_2)
{\tilde{f}}^2_{\Delta_{r_3-1}}(x_1,x_2){\tilde{f}}^1_{\Delta_{r_2-1}}(x_2), \\
p_{(BD,ABCD)}(x_1,x_2) & = & x_1x_2{\tilde{f}}^2_{\Delta_{r_1-1}}(x_1,x_2)
{\tilde{f}}^2_{\Delta_{r_2-1}}(x_1,x_2){\tilde{f}}^2_{\Delta_{r_3-1}}(x_1,x_2), \\
p_{(CD,ABCD)}(x_1,x_2) & = & x_1x_2{\tilde{f}}^2_{\Delta_{r_1-1}}(x_1,x_2)
{\tilde{f}}^2_{\Delta_{r_2-1}}(x_1,x_2){\tilde{f}}^2_{\Delta_{r_3-1}}(x_1,x_2).
\end{eqnarray*}

\begin{figure}
\centering
\includegraphics[scale=.7]{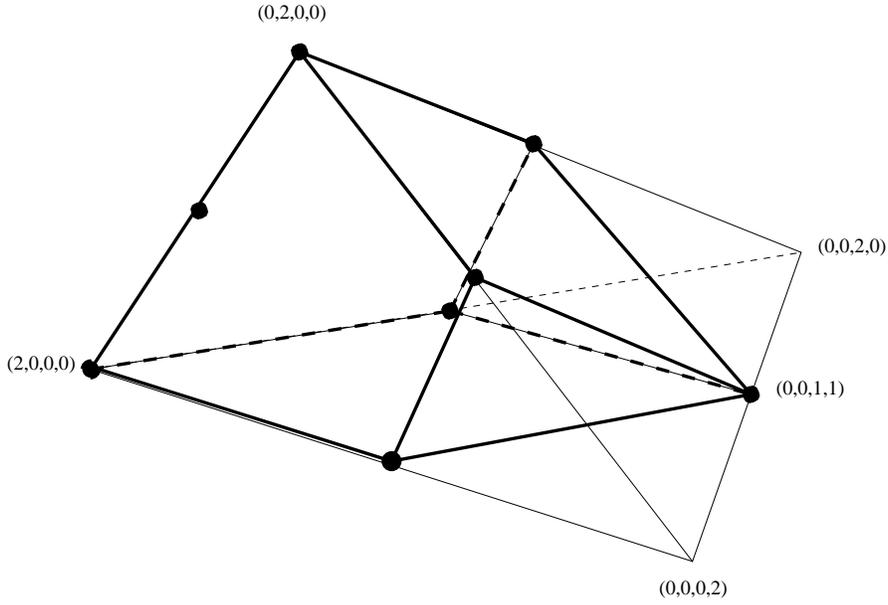}
\caption{$\Delta_{F_1} + \Delta_{F_2}$ where $F_1 = \{1,2,3\}$ and $F_2 = \{1,2,4\}$.}
\label{fig:contracted-cube}
\end{figure}

Finally, for $X\in\{A,B,C,D\}$, and $\tilde{A} = (A_1,A_2) = (X,ABCD)$
we have 
\[
\partial(\dim(\tilde{A})) = (\dim(X),\dim(ABCD)-\dim(X)) = (0,2),
\]
$\supp(A_1) = \supp(X)$, $\supp(A_2)\setminus\supp(A_1) = \{1,2,3\}\setminus\supp(X)$,
and hence
\[
p_{\tilde{A}}(x_1,x_2) = p_{(X,ABCD)}(x_1,x_2) = 
x_2^2\prod_{{\alpha}_1 \in\supp(X)}{\tilde{f}}^2_{\Delta_{r_{{\alpha}_1}-1}}(x_1,x_2)
\prod_{{\alpha}_2 \in\{1,2,3\}\setminus\supp(X)}{\tilde{f}}^1_{\Delta_{r_{{\alpha}_2}-1}}(x_2).
\]
From this we obtain
\begin{eqnarray*}
p_{(A,ABCD)}(x_1,x_2) & = & x_2^2{\tilde{f}}^2_{\Delta_{r_1-1}}(x_1,x_2)
{\tilde{f}}^1_{\Delta_{r_2-1}}(x_2){\tilde{f}}^1_{\Delta_{r_3-1}}(x_2), \\
p_{(B,ABCD)}(x_1,x_2) & = & x_2^2{\tilde{f}}^2_{\Delta_{r_1-1}}(x_1,x_2)
{\tilde{f}}^2_{\Delta_{r_2-1}}(x_1,x_2){\tilde{f}}^1_{\Delta_{r_3-1}}(x_2), \\
p_{(C,ABCD)}(x_1,x_2) & = & x_2^2{\tilde{f}}^2_{\Delta_{r_1-1}}(x_1,x_2)
{\tilde{f}}^2_{\Delta_{r_3-1}}(x_1,x_2){\tilde{f}}^1_{\Delta_{r_2-1}}(x_2), \\
p_{(D,ABCD)}(x_1,x_2) & = & x_2^2{\tilde{f}}^2_{\Delta_{r_2-1}}(x_1,x_2)
{\tilde{f}}^2_{\Delta_{r_3-1}}(x_1,x_2){\tilde{f}}^1_{\Delta_{r_1-1}}(x_2).
\end{eqnarray*}
Hence, collecting all the data for $P_{\cal{F}} = \Delta_{F_1} + \Delta_{F_2}$ with
$r_1 = |F_1\cap F_2|$, $r_2 = |F_1\setminus F_2|$ and $r_3 = |F_2\setminus F_1|$
we have by Theorem~\ref{thm:flag-main} that 
\[
{\tilde{f}}^2_{P_{\cal{F}}}(x_1,x_2) = 
\sum_{\tilde{A}\in {\mathbf{C}}_2(P(2))}p_{\tilde{A}}(x_1,x_2),
\]
where the sum is taken over the 25 different 2-chains $\tilde{A}$ given
in (\ref{eqn:2-chain}). Substituting we obtain the following.
\begin{corollary}
\label{cor:2-chain}
If $F_1,F_2\subseteq [r]$ and
$r_1 = |F_1\cap F_2|$, $r_2 = |F_1\setminus F_2|$ and $r_3 = |F_2\setminus F_1|$
then for $P_{\cal{F}} = \Delta_{F_1} + \Delta_{F_2}$ We have 
\begin{eqnarray*}
{\tilde{f}}^2_{P_{\cal{F}}}(x_1,x_2) 
& = &   (2x_1 + x_1^2 + 2x_1x_2){\tilde{f}}_1^2{\tilde{f}}_2^2{\tilde{f}}_3^2   
      + (x_2 + x_1x_2 + x_2^2)({\tilde{f}}_1^2{\tilde{f}}_2^2{\tilde{f}}_3^1 
+{\tilde{f}}_1^2{\tilde{f}}_3^2{\tilde{f}}_2^1) \\   
&   & + (2x_2 + x_2^2){\tilde{f}}_2^2{\tilde{f}}_3^2{\tilde{f}}_1^1 
      + x_2^2{\tilde{f}}_1^2{\tilde{f}}_2^1{\tilde{f}}_3^1 \\
&   & + (1 + x_1 + x_2)({\tilde{f}}_1^2{\tilde{f}}_2^2 + {\tilde{f}}_1^2{\tilde{f}}_3^2)
+ x_2({\tilde{f}}_1^2{\tilde{f}}_2^1 + {\tilde{f}}_1^2{\tilde{f}}_3^1)
+ {\tilde{f}}_2^2{\tilde{f}}_3^2 
+ {\tilde{f}}_1^2,
\end{eqnarray*}
where ${\tilde{f}}_i^2 = {\tilde{f}}^2_{\Delta_{r_i-1}}(x_1,x_2)$ 
and ${\tilde{f}}_i^1 = {\tilde{f}}^1_{\Delta_{r_i-1}}(x_2)$ for $i\in\{1,2,3\}$
are as in Observation~\ref{obs:flag-poly-simp}.
\end{corollary}
\noindent{\sc Example:}
For $F_1 = \{1,2,3\}$ and $F_2 = \{1,2,4\}$ we have
$r_1 = 2$ and $r_2 = r_3 = 1$ and we get in this case that
\[
{\tilde{f}}^2_{\Delta_{F_1} + \Delta_{F_2}}(x_1,x_2) = 
7 + 11x_1 + 6x_1^2 + x_1^3 +
22x_2 + 22x_1x_2 + 6x_1^2x_2 +
22x_2^2 + 11x_1x_2^2 +
7x_2^3,
\]
\comment{and hence
\[
{\tilde{\phi}}^2_{\Delta_{F_1} + \Delta_{F_2}}(y_1,y_2) = 
7 + 11y_1y_2 + 6y_1^2y_2^2 + y_1^3y_2^3 +
22y_2 + 22y_1y_2^2 + 6y_1^2y_2^3 +
22y_2^2 + 11y_1y_2^3 +
7y_2^3,
\]}
something that can be seen by inspection, since the Minkowski sum $\Delta_{F_1} + \Delta_{F_2}$
in this case has the same combinatorial structure as a cube where one side/edge
has been contracted to a vertex. 
The polytope $\Delta_{F_1} + \Delta_{F_2}$
has been drawn in Figure~\ref{fig:contracted-cube}.  

\vspace{3 mm}

\noindent{\sc Example:}
Using the setup from Corollary~\ref{cor:2-chain}
above, we note that the total number of vertices
of $P_{\cal{F}} = \Delta_{F_1} + \Delta_{F_2}$ is the
number of $(0,0)$ 2-chains of faces, which is precisely
the constant term in ${\tilde{f}}^2_{P_{\cal{F}}}(x_1,x_2)$
given by ${\tilde{f}}^2_{P_{\cal{F}}}(0,0)$.
Since ${\tilde{f}}^2_{\Delta_{r_i-1}}(x_1,x_2) = ((x_1+x_2+1)^{r_i} - (x_2+1)^{r_i})/x_1$ 
and ${\tilde{f}}^1_{\Delta_{r_i-1}}(x_2) = ((x_2+1)^{r_i} - 1)/x_2$
by Observation~\ref{obs:flag-poly-simp}, we get 
${\tilde{f}}^2_{\Delta_{r_i-1}}(0,0) = {\tilde{f}}^1_{\Delta_{r_i-1}}(0) = r_i$
for each $i$, and hence by Corollary~\ref{cor:2-chain} the number of vertices is 
in this case
\[
{\tilde{f}}^2_{P_{\cal{F}}}(0,0) = r_1r_2 + r_1r_3 + r_2r_3 + r_1.
\]
We will revisit this example in Section~\ref{sec:f-poly-2}.

\section{Generalizing the free join of polytopes}
\label{sec:join}

In this section we describe how viewing the free join of two polytopes
from a different angle, naturally generalizes to an action of a
polytope on a set of polytopes. We then present the $\ell$-flag polynomial
for the resulting polytope and finally apply that to the
Minkowski sum of simplices.

For polytopes $P$ and $Q$ recall that the {\em free join} 
(or just the {\em join} as it is sometimes called when there
is no danger of ambiguity~\cite {Ziegler})  $P*Q$
is a polytope that is obtained by positioning $P$ and $Q$ into
skew affine subspaces and then take the convex hull $\conv(P,Q)$.
The faces of $P*Q$ are given by (i) type $P$: the faces of $P$,
(ii) type $Q$: the faces of $Q$ and (iii) type $PQ$: the disjoint union $A\cup B$ 
where $A$ is a face of $P$ and $B$ a face of $Q$. The dimension
of a face of type $P$ or $Q$ is equal to its dimension
as a face in $P$ or $Q$, but the dimension of a face $A\cup B$ 
of type $PQ$ is equal to $\dim(A) + \dim(B) + 1$ because
of the convex combination of points of $A$ and points of $B$. 
The three types of faces of $P*Q$ is a partition of its faces.  
The free join of a polytope
$P$ and a singleton $\tilde{x}_0$ is called the {\em pyramid over $P$}
and denoted $\Pyr(P)$~\cite{Ziegler}. 

For a chain $C_1\subseteq C_2 \subseteq\cdots\subseteq C_{\ell}$
of faces of $P*Q$ there are five types of chains depending on
the type of the starting and terminal face:
types $(P,P)$, $(P,PQ)$, $(PQ,PQ)$, $(Q,PQ)$ and $(Q,Q)$.
The $\ell$-flag polynomial for type $(P,P)$ and $(Q,Q)$
chains is ${\tilde{f}}^{\ell}_{P}(\tilde{x})$ and 
${\tilde{f}}^{\ell}_{Q}(\tilde{x})$ respectively. 
The $\ell$-flag polynomial for a type $(PQ,PQ)$ chain
is the same as for $P\times Q$ except for the fact
that the dimension of each face $A\cup B$ is one more
that that of $A\times B$. Hence, the $\ell$-flag polynomial
must be shifted upward by one at each level of the chain.
By definition of the flag polynomial, that is obtained
by multiplying by the the first variable $x_1$ and we obtain
$x_1{\tilde{f}}^{\ell}_{P}(\tilde{x}){\tilde{f}}^{\ell}_{Q}(\tilde{x})$.
For a type $(P,PQ)$ chain, there are $\ell - 1$ possible places
in the chain where we encounter the first type $PQ$ 
face, say at step $i\in \{2,\ldots,\ell\}$. Each such
chain is determined by an $\ell$-chain of faces in $P$
(the $A_i$s) and an additional $(\ell - i + 1)$-chain in $Q$ (the $B_i$s) 
and has the form
\begin{equation}
\label{eqn:(P,PQ)}
A_1\subseteq\cdots\subseteq A_{i-1}\subseteq A_i*B_i\subseteq\cdots\subseteq A_{\ell}*B_{\ell}.
\end{equation}
Again, since the dimension of each $A_j*B_j$ is one more than that of $A_j\times B_j$,
we must shift the dimension upward by one from and including at level $i$,
which is obtained by multiplying by $x_i$. Hence, the $\ell$-flag polynomial 
for all chains of the form (\ref{eqn:(P,PQ)}) is given by 
$x_i{\tilde{f}}^{\ell}_{P}(\tilde{x}){\tilde{f}}^{\ell-i+1}_{Q}(x_i,\ldots,x_{\ell})$.
Since the case for type $(Q,PQ)$ chains is treated in the same way we
obtain the following.
\begin{proposition}
\label{prp:P*Q}
For polytopes $P$ and $Q$ and a positive integer $\ell\in \nats$, the
$\ell$-flag polynomial for the free join $P*Q$ 
is given by
\begin{eqnarray*}
{\tilde{f}}^{\ell}_{P*Q}(\tilde{x}) & = & {\tilde{f}}^{\ell}_{P}(\tilde{x})
   \left(1 + \sum_{i=2}^{\ell}x_i{\tilde{f}}^{\ell-i+1}_{Q}(x_i,\ldots,x_{\ell})\right) \\
                                  & + & {\tilde{f}}^{\ell}_{Q}(\tilde{x})
   \left(1 + \sum_{i=2}^{\ell}x_i{\tilde{f}}^{\ell-i+1}_{P}(x_i,\ldots,x_{\ell})\right) \\
                                  & + & x_1{\tilde{f}}^{\ell}_{P}(\tilde{x}){\tilde{f}}^{\ell}_{Q}(\tilde{x}).
\end{eqnarray*}
\end{proposition}
Note that for the empty polytope $\emptyset$ we have $P*\emptyset = P$ and $\emptyset*Q = Q$.
Also, for each $\ell$ we have ${\tilde{f}}^{\ell}_{\emptyset}(\tilde{x}) = 0$.

When $Q$ is a singleton $Q = \Delta_0$, then  ${\tilde{f}}^{\ell}_{Q}(\tilde{x}) = 1$ 
for each $\ell$, and hence we obtain the following.
\begin{corollary}
\label{cor:flag-poly-pyr}
For a polytope $P$ and $\ell\in\nats$, the $\ell$-flag polynomial 
for the pyramid $\Pyr(P)$ over $P$ is given by
\[
{\tilde{f}}^{\ell}_{\Pyr(P)}(\tilde{x}) = 
{\tilde{f}}^{\ell}_{P}(\tilde{x})(1 + x_1 + \cdots + x_{\ell}) + 1 +
\sum_{i=2}^{\ell}x_i{\tilde{f}}^{\ell-i+1}_{P}(x_i,\ldots,x_{\ell}).
\]
\end{corollary}

\noindent{\sc Remark:} For a $d$-polytope the $d$-flag-polynomial contains all the information
of the flag vector and the $\fatc\fatd$-index of $P$ with many redundancies
(i.e.~the coefficient of different monomials of the flag polynomial will count 
the same number of strict chains). Working with the most compact encoding
of the flag vector, elegant formulae for computing the $\fatc\fatd$-indices 
of the prism $\Pri(P) = P\times\Delta_1$ and the pyramid $\Pyr(P) = P*\Delta_0$
in terms of the $\fatc\fatd$-index $\Psi(P)$ of $P$ and some linear derivations
${\reals}\langle \fatc,\fatd\rangle \rightarrow {\reals}\langle \fatc,\fatd\rangle$
are given in~\cite{Ehrenborg-Readdy}. But even with these nice recursions,
it is hard to obtain explicit formulae for the $\fatc\fatd$-index of specific
polytopes, like the simplex $\Delta_d = \Pyr(\Delta_{d-1})$, largely since
derivations on non-commuting monomials in ${\reals}\langle \fatc,\fatd\rangle$ 
expand exponentially in repeated steps and do not simplify easily. These recursions 
are generalized to a recursion to compute the $\fatc\fatd$-indices for the 
free join and the Cartesian product in~\cite{Ehrenborg-Fox}.

\vspace{3 mm}

We now generalize the idea of the free join and investigate the effect
on the corresponding flag polynomial from Proposition~\ref{prp:P*Q}. 
First note that a concrete geometric realization of the free join $Q_1*Q_2$
can be obtained by an action $(Q_1,Q_2) \mapsto \Delta_1(Q_1,Q_2)$
of the 1-simplex $\Delta_1$ on the polytopes $Q_1$ and $Q_2$.
Here $\Delta_1$ is the standard
geometric 1-dimensional simplex $\Delta_{[2]} = \{(x_1,x_2) : x_1, x_2 \geq 0, x_1 + x_2 = 1\}$,
and we attach 
$Q_1$ to one endvertex $(1,0)$ and $Q_2$ to the other $(0,1)$. Here we also choose an embedding
of each $Q_i$ in $\{\tilde{x}_i \in {\reals}^{d_i+1} : \sum x_{i\/j} = 1\}$ 
where $d_i = \dim(Q_i)$, and the action is given by
\[
\Delta_1(Q_1,Q_2) := \{(t_1\tilde{x}_1,t_2\tilde{x}_2) \in {\reals}^{d_1 + d_2 + 2} : 
  t_1, t_2 \geq 0, \ t_1 + t_2 = 1,\ \tilde{x}_1\in Q_1,\ \tilde{x}_2\in Q_2\}.
\]
Note that for a nonempty polytope $Q$ we have 
$\Delta_1(Q,\emptyset) = \Delta_1(\emptyset,Q) = Q$. However,
if both $Q_1$ and $Q_2$ are nonempty, then 
every face of $\Delta_1(Q_1,Q_2)$ has the form $A(F_1,F_2)$ where
each $F_i$ is a face of $Q_i$ or $Q_i$ itself, and $A$ is a face of $\Delta_1$
or $\Delta_1$ itself. 

Each $\ell$-chain of faces of $\Delta_1(Q_1,Q_2)$ yields,
by projection,
an $\ell$-chain of faces $A_1\subseteq\cdots\subseteq A_{\ell}$ of $\Delta_1$, and hence
we have a partition of the $\ell$-chains of faces of $\Delta_1(Q_1,Q_2)$ where each part
consists of the chains that yield the same chain of $\Delta_1$. Note that each proper 
face $A$ of $\Delta_1$ has $\supp(A) \in \{ \{1\}, \{2\}, \{1,2\}\}$ and hence
each $\ell$-chain $\tilde{A}$ of faces of $\Delta_1$ has 
${\tilde{x}}^{\partial(\dim(\tilde{A}))} \in \{1,x_1,\ldots,x_{\ell}\}$.
Writing $Q_1$ for $P$ and $Q_2$ for $Q$ in Proposition~\ref{prp:P*Q},
the $\ell$-flag polynomial in Proposition~\ref{prp:P*Q}
for $Q_1*Q_2$ equals that of is geometric realization $\Delta_1(Q_1,Q_2)$, and we obtain
\begin{equation}
\label{eqn:form2}
{\tilde{f}}^{\ell}_{Q_1*Q_2}(x_1,\ldots,x_{\ell}) = 
\sum_{\tilde{A}\in {\mathbf{C}}_{\ell}(\Delta_1)}
\left(
{\tilde{x}}^{\partial(\dim(\tilde{A}))}
\prod_{i=1}^{\ell}
\left(
\prod_{{\alpha}_i \in\supp(A_i)\setminus \supp(A_{i-1})}
{\tilde{f}}^{\ell-i+1}_{Q_{{\alpha}_i}}(x_i,\ldots,x_{\ell})
\right)
\right).
\end{equation}

\noindent{\sc Conventions:} (i) Call an embedding of a polytope 
$P$ in ${\reals}^{\dim(P)+1}$ {\em simplicial} if
$P \subseteq \Delta_{\dim(P)} = \Delta_{[\dim(P)+1]}$.
(ii) By a {\em full face} of a polytope we will
mean either a face or the entire polytope itself. 

\vspace{3 mm}

The repeated free $m$-join $Q_1*\cdots *Q_m$ of polytopes $Q_i$s is well defined since
the operation $*$ is associative and commutative, and it has the empty polytope $\emptyset$
as the identity element. This operation $*$ makes the collection $\cal{Q}$ 
of all polytopes $Q$ a commutative semigroup $(\cal{Q},*)$. 
\begin{observation}
\label{obs:m-fold}
Let $Q_1,\ldots, Q_m$ be polytopes (the empty polytope $\emptyset$ is allowed).  
The $m$-fold free join $Q_1*\cdots*Q_m$ is a polytope 
of dimension $m - 1 + \sum_{i=1}^m \dim(Q_i)$ in which each full face 
has the the form $F_1*\cdots*F_m$ where each $F_i$
is a full face of $Q_i$ or the empty face. As $\dim(\emptyset) = -1$, 
the dimension of each face is given by
$\dim(F_1*\cdots*F_m) = m - 1 + \sum_{i=1}^m \dim(F_i)$.
\end{observation}
One nice advantage of (\ref{eqn:form2}) over the formula in 
Proposition~\ref{prp:P*Q} is that it generalizes
naturally in light of Observation~\ref{obs:m-fold}. This is
since a geometric realization of 
the $m$-fold free join $Q_1*\cdots*Q_m$ can, as above, be obtained
by an action $(Q_1,\ldots,Q_m)\mapsto \Delta_{m-1}(Q_1,\ldots,Q_m)$
of the standard $(m-1)$-simplex $\Delta_{m-1}$ on $Q_1,\ldots, Q_m$.
We attach each $Q_i$ to an endvertex $\tilde{e}_i$ of $\Delta_{m-1}$
and so 
\[
\Delta_{m-1}(Q_1,\ldots,Q_m) := \{(t_1\tilde{x}_1,\ldots,t_m\tilde{x}_m) \in {\reals}^{m + \sum d_i} : 
\sum_i t_i = 1,\   
t_i\geq 0,\ 
\tilde{x}_i\in Q_i \mbox{ for each }i\}.
\]
Here we do assume each $Q_i$ to be a nonempty polytope, and so
each full face of $\Delta_{m-1}(Q_1,\ldots,Q_m)$ has the form 
$A(F_1,\ldots,F_m)$ where $A$ is a full face of $\Delta_{m-1}$
and each $F_i$ is a full face of $Q_i$. 

As $\Delta_{m-1}(Q_1,\ldots,Q_m)$ is the geometric realization of
$Q_1*\cdots*Q_m$, their flag-$\ell$ polynomials are identical, and
we have the following corollary. 
\begin{corollary}
\label{cor:flag-m-fold}
For $\ell\in\nats$ and simplicially embedded polytopes $Q_1,\ldots,Q_m$, the flag-$\ell$ polynomial
for the $m$-fold free join $Q_1*\cdots*Q_m$ is given by
\[
{\tilde{f}}^{\ell}_{Q_1*\cdots *Q_m}(x_1,\ldots,x_{\ell}) = 
\sum_{\tilde{A}\in {\mathbf{C}}_{\ell}(\Delta_{m-1})}
\left(
{\tilde{x}}^{\partial(\dim(\tilde{A}))}
\prod_{i=1}^{\ell}
\left(
\prod_{{\alpha}_i \in\supp(A_i)\setminus \supp(A_{i-1})}
{\tilde{f}}^{\ell-i+1}_{Q_{{\alpha}_i}}(x_i,\ldots,x_{\ell})
\right)
\right).
\]
\end{corollary}

We now consider a generalization $P(Q_1,\ldots,Q_m)$
of $\Delta_{m-1}(Q_1,\ldots,Q_m)$ where $P$ is now an arbitrary
geometric polytope of dimension $m-1$ with $\supp(P) = \{1,\ldots,m\}$. Here
we assume $P = \conv(\tilde{p}_1,\ldots,\tilde{p}_n)$ to be
simplicially embedded in ${\reals}^m$ and $\tilde{p}_i = (p_{i\/1},\ldots,p_{i\/m})$
for each $i$. The polytope $P$ acts on simplicially embedded $Q_1,\ldots,Q_m$ as 
$(Q_1,\ldots,Q_m)\mapsto P(Q_1,\ldots,Q_m)$ and
the action is given by the following. 
\begin{definition}
\label{def:P-acts}
For a simplicially embedded polytope $P = \conv(\tilde{p}_1,\ldots,\tilde{p}_n)$ 
with $\supp(P) = \{1,\ldots,m\}$ and simplicially embedded
polytopes $Q_1,\ldots,Q_m$ the action of $P$ is given by
\begin{eqnarray*}
P(Q_1\ldots,Q_m) & := & \{ t_1(p_{1\/1}\tilde{x}_1,\ldots,p_{1\/m}\tilde{x}_m) + \cdots +
                          t_n(p_{n\/1}\tilde{x}_1,\ldots,p_{n\/m}\tilde{x}_m) : \\
                 &   &  0\leq t_i \leq 1, \ \sum_i t_i = 1, \ \tilde{x}_i \in Q_i \mbox{ for each }i\}.
\end{eqnarray*}
\end{definition}
Similar to Observation~\ref{obs:m-fold} we have the following
generalization.
\begin{proposition}
\label{prp:P-acts}
The dimension of the polytope $P(Q_1,\ldots,Q_m)$ from Definition~\ref{def:P-acts}
is given by $m - 1 +  \sum_{i=1}^m \dim(Q_i)$ and it is 
simplicially embedded. Each full 
face has the form $A(F_1,\ldots,F_m)$ where $A$ is a full face of $P$
and each $F_i$ is a full face of $Q_i$.
Also 
\[
\dim(A(F_1,\ldots,F_m)) = \dim(A) + \sum_{i\in\supp(A)}\dim(F_i).
\]
\end{proposition}
\begin{proof}
For each $i$ let $d_i = \dim(Q_i)$ and $\tilde{x}_i \in{\reals}^{d_i+1}$
a tuple of variables and $\tilde{x}^* = (\tilde{x}_1,\ldots,\tilde{x}_m) \in {\reals}^{m+\sum_id_i}$. 
Each full face of $P(Q_1\ldots,Q_m)$ is determined
as the maximum of a functional $\tilde{x}^* \mapsto \tilde{a}^*\cdot\tilde{x}^*$.
For a given $\tilde{a}^* = (\tilde{a}_1,\ldots,\tilde{a}_m)$ and $i$,
let $F_i$ be the face of $Q_i$ that maximizes 
$\tilde{x}_i \mapsto \tilde{a}_i\cdot\tilde{x}_i$ and assume
this maximum is $M_i$. Let $A$ be the face of $P$ 
that is maximizes $(y_1,\ldots,y_m)\mapsto \sum_iM_iy_i$.
Then the face that maximizes $\tilde{x}^* \mapsto \tilde{a}^*\cdot\tilde{x}^*$
is $A(F_1,\ldots,F_m)$. 
\end{proof}
\noindent{\sc Remark:} Note that if $i\not\in\supp(A)$ then
$A(F_1,\ldots,F_m) = A(F_1,\ldots,F_{i-1},Q_i,F_{i+1},\ldots,F_m)$
regardless of the face $F_i$.

\vspace{3 mm}

We now briefly discuss how we interpret $P(Q_1\ldots,Q_m)$ when 
$Q_i=\emptyset$ is allowed: let $i\in\{1,\ldots,m\}$ and assume that 
$Q_i = \emptyset$ is the empty polytope.
As $P = \conv(\tilde{p}_1,\ldots,\tilde{p}_n)$ is simplicially embedded in ${\reals}^m$,
we can assume that $i\not\in\supp(\tilde{p}_1,\ldots,\tilde{p}_h)$
and $i\in\supp(\tilde{p}_{h+1},\ldots,\tilde{p}_n)$.
In this case let $P_{\hat{\imath}} = \conv(\pi_{\hat{\imath}}(\tilde{p}_1),\ldots,\pi_{\hat{\imath}}(\tilde{p}_h))$
where $\pi_{\hat{\imath}} : {\reals}^m \rightarrow {\reals}^{m-1}$ is the projection
obtained by omitting the $i$-th coordinate. In this case we make the
following definition.
\begin{equation}
\label{eqn:Q=empty}
P(Q_1\ldots,Q_m) = P(Q_1\ldots,Q_{i-1},\emptyset,Q_{i+1},\ldots,Q_m)
:= P_{\hat{\imath}}(Q_1,\ldots,Q_{i-1},Q_{i+1},\ldots,Q_m).
\end{equation}
By (\ref{eqn:Q=empty}) we can therefore assume that each $Q_i$ in
$P(Q_1\ldots,Q_m)$ is a proper nonempty polytope.

\vspace{3 mm}

We now state and prove our main result in this section, that the $\ell$-flag 
polynomial for $P(Q_1,\ldots,Q_m)$
has the same form as displayed in Corollary~\ref{cor:flag-m-fold}.
\begin{theorem}
\label{thm:P-acts-flag}
For $\ell\in\nats$ and simplicially embedded polytopes $Q_1,\ldots,Q_m$, the flag-$\ell$ polynomial
of the polytope $P(Q_1,\ldots,Q_m)$ from Definition~\ref{def:P-acts}
is given by
\[
{\tilde{f}}^{\ell}_{P(Q_1,\ldots,Q_m)}(x_1,\ldots,x_{\ell}) = 
\sum_{\tilde{A}\in {\mathbf{C}}_{\ell}(P)}
\left(
{\tilde{x}}^{\partial(\dim(\tilde{A}))}
\prod_{i=1}^{\ell}
\left(
\prod_{{\alpha}_i \in\supp(A_i)\setminus \supp(A_{i-1})}
{\tilde{f}}^{\ell-i+1}_{Q_{{\alpha}_i}}(x_i,\ldots,x_{\ell})
\right)
\right).
\]
\end{theorem}
We will prove Theorem~\ref{thm:P-acts-flag} in a few steps and
use that to prove Theorem~\ref{thm:flag-main}. For simplicity
we will write $\mathbf{X} = (X_1,\ldots, X_m)$. Assuming that
$\dim(Q_i) = d_i$ for each $i$ we have linear maps
$h = h_{P(\mathbf{Q})} : {\reals}^{m+\sum d_i} \rightarrow {\reals}^m$ and
$\pi_{\tilde{\imath}} : {\reals}^{m+\sum d_i} \rightarrow {\reals}^{d_i+1}$
given by
\[
h(\tilde{x}_1,\ldots,\tilde{x}_m) = \left(\sum_{i=1}^{d_1}x_{1\/i},\ldots,\sum_{i=1}^{d_m}x_{m\/i}\right),
\]
and 
\[
\pi_{\tilde{\imath}}(\tilde{x}_1,\ldots,\tilde{x}_m) = \tilde{x}_i
\]
respectively.
\begin{lemma}
\label{lmm:AA'}
For faces of the polytope $P(\mathbf{Q})$ we have
$A(\mathbf{F})\subseteq A'(\mathbf{F}')$ if and only if
$A\subseteq A'$ are subfaces in $P$ and $F_i\subseteq F_i'$ are subfaces
in $Q_i$ for each $i\in\supp(A)$.
\end{lemma}
\begin{proof}
We will prove the more involved direction. Assume that
$A(\mathbf{F})\subseteq A'(\mathbf{F}')$. Since each of the $Q_i$
are simplicially embedded, we have 
$A = h(A(\mathbf{F}))\subseteq  h(A'(\mathbf{F}')) = A'$. 

For $i\in\supp(A)$, then by the definition of 
$P(\mathbf{Q}) = P(Q_1,\ldots,Q_m)$, we have $\pi_{\tilde{\imath}}(A(\mathbf{F})) = [a_i,b_i]F_i$
where $\pi_i(A) = [a_i,b_i]$ and $0\leq a_i < b_i\leq 1$. (Here $\pi_i$ is the projection
onto the $i$-th coordinate). Since we have that $A\subseteq A'$, then for each $i$
we have $[a_i,b_i]\subseteq [a_i',b_i']$ or $a_i'\leq a_i<b_i\leq b_i'$,
or 
\[
[a_i,b_i]F_i = \pi_{\tilde{\imath}}(A(\mathbf{F}))
     \subseteq \pi_{\tilde{\imath}}(A'(\mathbf{F}')) = [a_i',b_i']F_i'.
\]
For $c_i\in [a_i,b_i]$ we have $c_iF_i\subseteq [a_i,b_i]F_i\subseteq [a_i',b_i']F_i'$.
Since $c_iF_i = \{\tilde{x}\in [a_i,b_i]F_i : \sum_jx_j = c_i\}$, then
$c_iF_i\subseteq [a_i',b_i']F_i'\cap\{\tilde{x} : \sum_jx_j = c_i\} = c_iF_i'$
(Here $Q_i$ is simplicially embedded and so are all of its faces).
Therefore 
\[
F_i = \frac{1}{c_i}(c_iF_i)\subseteq \frac{1}{c_i}(c_iF_i') = F_i',
\]
which completes the proof of the lemma.
\end{proof}
By Proposition~\ref{prp:P-acts} each $\ell$-chain of faces of $P(\mathbf{Q})$
has the form
\begin{equation}
\label{eqn:l-chain}
A_1(\mathbf{F}_1)\subseteq A_2(\mathbf{F}_2)\subseteq\cdots\subseteq A_{\ell}(\mathbf{F}_{\ell})
\end{equation}
which by Lemma~\ref{lmm:AA'} firstly yields a unique $\ell$-chain 
$\tilde{A} : A_1\subseteq A_2\subseteq\cdots\subseteq A_{\ell}$ of faces of $P$,
and hence $\supp(A_1)\subseteq \supp(A_2)\subseteq \cdots\subseteq \supp(A_{\ell})$.
By a suitable permutation of the $\{1,\ldots,m\}$ we may assume that 
$\supp(A_1) < \supp(A_2)\setminus\supp(A_1) < \cdots < \supp(A_{\ell})\setminus\supp(A_{\ell -1})$
(where $X < Y$ means that $x < y$ for all $x\in X$ and $y\in Y$).
Secondly, by this mentioned assumption on the order on the $\supp(A_i)$, then
for each fixed $\ell$-chain $\tilde{A}\in {\mathbf{C}}_{\ell}(P)$
(\ref{eqn:l-chain}) yields a Young tableau (or Ferres diagram) of boxes
in $[\ell]\times[m]$ such that for each $i\in [\ell]$ the boxes 
numbered $(i,1),\ldots, (i,a_i)$, where $a_i = |\supp(A_i)|$, 
contain the faces $F_{i\/1},F_{i\/2},\ldots,F_{i\/a_i}$
from $\mathbf{F}_i$ that are indexed by numbers of $\supp(A_i)$.
Let $1=i_1 < i_2 < \cdots < i_{\lambda}\leq \ell$ be the complete
list of numbers from $[\ell]$ where $a_{i_j-1} < a_{i_j}$ for each 
$j=1,\ldots,\lambda$. Hence, for each such $j$ we obtain $a_{i_j}-a_{i_j-1}$
chains of faces 
\[
F_{i_j\/{\alpha}} \subseteq F_{i_j+1\/{\alpha}}\subseteq \cdots \subseteq F_{\ell\/{\alpha}},
\]
in $Q_{\alpha}$, one for each $\alpha\in\supp(A_{i_j})\setminus\supp(A_{i_j-1})$.
Conversely, every $\ell$-chain as in (\ref{eqn:l-chain}) is presented
by such a Young tableau in which the horizontal sequence of boxes at hight $\alpha$
yields a chain of faces in $Q_{\alpha}$. This together with 
Proposition~\ref{prp:P-acts} yields the following lemma.
\begin{lemma}
\label{lmm:key-chain}
For $\ell\in\nats$ let $\tilde{A}\in {\mathbf{C}}_{\ell}(P)$ be a fixed
$\ell$-chain of faces of $P$. Then the coefficient of ${\tilde{x}}^{\partial{\tilde{s}}}$
in the polynomial given by
\begin{equation}
\label{eqn:key-poly}
{\tilde{x}}^{\partial(\dim(\tilde{A}))}
\prod_{j=1}^{\lambda}
\left(
\prod_{\alpha_j\in\supp(A_{i_{j}})\setminus \supp(A_{i_{j-1}})}
{\tilde{f}}^{\ell-i_j+1}_{Q_{\alpha_j}}(x_{i_j},\ldots,x_{\ell})
\right)
\end{equation}
is the number of $\ell$-chains $C_1\subseteq \cdots \subseteq C_{\ell}$ of faces 
of $P(Q_1,\ldots,Q_m)$ with $\dim(C_i) = s_i$ and $h(C_i) = A_i$.
\end{lemma}
With the convention that a product with an empty index set equals
one, (\ref{eqn:key-poly}) can be written as 
\[
{\tilde{x}}^{\partial(\dim(\tilde{A}))}
\prod_{i=1}^{\ell}
\left(
\prod_{{\alpha}_i \in\supp(A_i)\setminus \supp(A_{i-1})}
{\tilde{f}}^{\ell-i+1}_{Q_{{\alpha}_i}}(x_i,\ldots,x_{\ell})
\right).
\]
Theorem~\ref{thm:P-acts-flag} now follows directly from Lemma~\ref{lmm:key-chain}.

Now let $k,r\in\nats$ and consider a Minkowski sum $P_{\cal{F}}$ of simplices as in 
(\ref{eqn:Minkssum}). For $i\in [2^k-1]$ let $r_i$ be as before, the number
of elements of $R_i = \{j\in [r]: N_{\cal{F}}(j) = N_{{\cal{H}}(k)}(i)\}$
and so $r = r_1 + \cdots + r_{2^k-1}$. 
If $n\in\nats$ and $C$ is a convex set, then 
$C + \cdots + C \mbox{ ($n$ times) } = nC$. Hence, from this and by
Definitions~\ref{def:master} of the master polytope $P(k)$
and~\ref{def:P-acts} we have that
\begin{eqnarray}
P_{\cal{F}} & = & \sum_{F\in{\cal{F}}}\Delta_F \nonumber \\ 
            & = & \Delta_{F_1} + \cdots + \Delta_{F_k} \nonumber \\ 
            & = & P(k)(\Delta_{[r_1]},\Delta_{[r_2]},\ldots,\Delta_{[r_{2^k-1}]}) \nonumber \\ 
            & = & P(k)(\Delta_{r_1-1},\Delta_{r_2-1},\ldots,\Delta_{r_{2^k-1}-1}). \label{eqn:Minks=P(k)act}
\end{eqnarray}
Theorem~\ref{thm:flag-main} now follows from (\ref{eqn:Minks=P(k)act}) and
Theorem~\ref{thm:P-acts-flag}.

\section{Extremal considerations of the Minkowski sum of two simplices}
\label{sec:f-poly-2}

In this section we consider two extremal examples of the Minkowski sum
of two simplices, where we apply the flag polynomial. In our first 
example we obtain the $f$-polynomial 
of a general Minkowski sum of $k$ standard simplices in terms of the master 
polytope and use that to compute the
maximum number of $d$-faces of the Minkowski sum of two simplices when $k=2$. 
In our second example we compute the maximum 
discrepancy between the number of $(0,d)$-chains of faces of the Minkowski 
sum and the minimum number of such $(0,d)$-chains of a simple polytope 
of the same dimension and on the same number of vertices.

\subsection{First Example}

As before, $P_{\cal{F}}$ will be the Minkowski sum of $k$ simplices as in (\ref{eqn:Minkssum}).
Note that each face $Q$ of $P_{\cal{F}}$ can be viewed as a $1$-chain of faces.
Hence, letting $\ell = 1$ in Theorem~\ref{thm:flag-main} we get 
${\mathbf{C}}_{1}(P(k)) = \mathbf{f}(P(k))$, the set of faces of $P(k)$, and
we obtain a generating
function in terms of $x_1$ alone, which we denote by just $x$ in the
following corollary.
\begin{corollary}
\label{cor:main}
For a collection of $k$ ordered subsets ${\cal{F}} = (F_1,F_2,\ldots,F_k)$ of $[r]$,
the $f$-polynomial of the Minkowski sum $P_{\cal{F}}$ is given by
\[
f_{P_{\cal{F}}}(x) = 
\sum_{A\in \mathbf{f}(P(k))}\left(x^{\dim(A)}\prod_{i\in\supp(A)}f_{\Delta_{r_i-1}}(x)\right).
\]
\end{corollary}
We will now use Corollary~\ref{cor:main}
to investigate in more detail the Minkowski sum of two simplices, the $f$-polynomial
and then the maximum number of $d$-dimensional faces such a polytope can have.

By Corollary~\ref{cor:main} we have in particular for $k=2$ and ${\cal{F}} = (F_1,F_2)$ 
(and therefore $r_1 = |F_1\cap F_2|$, $r_2 = |F_1\setminus F_2|$ and $r_3 = |F_2\setminus F_1|$,) 
that the $f$-polynomial of $P_{\cal{F}} = \Delta_{F_1} + \Delta_{F_2}$ is given by 
\begin{eqnarray*}
f_{P_{\cal{F}}}(x) & = & (x^2 + 2x)f_{\Delta_{r_1-1}}(x)f_{\Delta_{r_2-1}}(x)f_{\Delta_{r_3-1}}(x) \\
                   & + & (x+1)\left( f_{\Delta_{r_1-1}}(x)f_{\Delta_{r_2-1}}(x) + f_{\Delta_{r_1-1}}(x)f_{\Delta_{r_3-1}}(x) \right)\\
                   & + & f_{\Delta_{r_2-1}}(x)f_{\Delta_{r_3-1}}(x) \\
                   & + & f_{\Delta_{r_1-1}}(x).
\end{eqnarray*}
Hence, together with Observation~\ref{obs:f-poly-standard} we have the following.
\begin{corollary}
\label{cor:k=2}
If $P_{\cal{F}} = \Delta_{F_1} + \Delta_{F_2}$ and $r = r_1+r_2+r_3$, 
where $r_1 = |F_1\cap F_2|$, $r_2 = |F_1\setminus F_2|$ and $r_3 = |F_2\setminus F_1|$,
then the $f$-polynomial of $P_{\cal{F}}$ is given by 
\[
f_{P_{\cal{F}}}(x) = \frac{1}{x^2}\left( (x+2)(x+1)^r - (x+1)^{r_1+r_2} - (x+1)^{r_1+r_3} - (x+1)^{r_2+r_3+1} + 1\right),
\]
and hence for each $d\in\{0,1,\ldots,r-1\}$ the number $f_d(P_{\cal{F}})$ of $d$-dimensional
faces of $P_{\cal{F}}$ is given by
\[
f_d(P_{\cal{F}}) = \binom{r+1}{d+2} + \binom{r}{d+2} - \binom{r_1+r_2}{d+2} - \binom{r_1+r_3}{d+2} - \binom{r_2+r_3+1}{d+2}.
\]
\end{corollary}
\noindent{\sc Remark:} Writing $f_{P_{\cal{F}}}(x)$ directly in terms of sets $F_1,F_2,F_3$ we obtain
\[
f_{P_{\cal{F}}}(x) = \frac{1}{x^2}\left( (x+2)(x+1)^{|F_1\cup F_2|} 
                   - (x+1)^{|F_1|} - (x+1)^{|F_2|} - (x+1)^{|F_1 \triangle F_2| + 1} + 1\right).
\]

For given $r\in\nats$ and $d\in\{0,1,\ldots,r-1\}$ we want to find the the maximum number
$f_d(r)$ of $d$-dimensional faces a Minkowski sum $P_{\cal{F}} = \Delta_{F_1} + \Delta_{F_2}$ of
two simplices with $F_1,F_2\subseteq \{1,\ldots, r\}$ can have. That is, we want to determine
\begin{equation}
\label{eqn:max-2d-faces}
f_{2;d}(r) = \max\{f_d(P_{\cal{F}}) : P_{\cal{F}} = \Delta_{F_1} + \Delta_{F_2}, \ \ r_1+r_2+r_3 = r\}.
\end{equation}
To do this, it is understood
that for a nonzero integer $c$ the map $\ints\rightarrow\ints$ defined as $x\mapsto \binom{x}{c}$
is given by $\binom{x}{c} = \frac{x(x-1)\cdots(x-c+1)}{c!}$ for $x\geq c$ and $\binom{x}{c} = 0$ 
for $x < c$.
\begin{lemma}
\label{lmm:binom-real}
For nonzero integers $c$ and $x<y$ we have that
\[
\binom{x}{c}+\binom{y}{c}\geq \binom{x+1}{c}+\binom{y-1}{c}.
\]
\end{lemma}
\begin{proof}
Considering each of the cases, $c\leq x<y$, $x<c<y$ and $x<y\leq c$, and
using the fact that $\binom{t+1}{c} -\binom{t}{c} = \binom{t}{c-1}$ when $t\geq c-1$,
the lemma follows.
\end{proof}
\begin{corollary}
\label{cor:min-123}
Let $K\in\nats$ and $x_1,x_2,x_3\geq 0$ be integers with $x_1+x_2+x_3 = K$
and 
\[
g(x_1,x_2,x_3) = \binom{x_1}{c}+\binom{x_2}{c}+\binom{x_3}{c}.
\]
\begin{itemize}
  \item If $K \geq 3c - 2$, then $g(x_1,x_2,x_3)$ will take it 
minimum value if, and only if, $|x_i-x_j|\leq 1$ for each $i,j$.  
  \item If $K \leq 3(c-1)$, then any $x_1,x_2,x_3\leq c-1$ summing
up to $K$ will yield the minimum value $0$ of $g(x_1,x_2,x_3)$.
In particular those with $|x_i-x_j|\leq 1$ for each $i,j$.
\end{itemize}
\end{corollary}
Clearly $f_d(P_{\cal{F}})$ from Corollary~\ref{cor:k=2} is maximized
iff 
\[
h(r_1,r_2,r_3) =  \binom{r_1+r_2}{d+2} + \binom{r_1+r_3}{d+2} + \binom{r_2+r_3+1}{d+2}.
\]
is minimized. Note that $(r_1+r_2) + (r_1+r_3) + (r_2+r_3+1) = 2r+1$. Considering 
each of the cases of $r \equiv 0,1,2 \pmod{3}$ we obtain the following
theorem, in which $[x]$ denotes the {\em nearest} integer to the real
number $x$. Note that for a fixed $x = a/b$, where $a,b$ are integers and $b$ is odd, then
$y = [x]$ is the unique integer that minimizes $|x-y|$.  
\begin{theorem}
\label{thm:max-2d-faces}
If $f_{2;d}(r)$ is the maximum number of $d$-faces in a Minkowski sum
$\Delta_{F_1} + \Delta_{F_2}$ of two simplices where $F_1, F_2\subseteq \{1,\ldots,r\}$
as given in (\ref{eqn:max-2d-faces}), then
\[
f_{2;d}(r) = 
\binom{r+1}{d+2} + \binom{r}{d+2} - \binom{\floor{(2r+1)/3}}{d+2} - \binom{[(2r+1)/3]}{d+2} - \binom{\ceil{(2r+1)/3}}{d+2}.
\]
\end{theorem}
\noindent{\sc Remarks:} (i) In the case when $2(r-1)/3\leq d\leq r-1$ we note that 
$f_{2;d}(r) = \binom{r+1}{d+2} + \binom{r}{d+2}$. Moreover, any 
Minkowski sum $P_{\cal{F}} = \Delta_{F_1} + \Delta_{F_2}$
such that $r_1+r_2+r_3 = r$, $r_1\geq r-d$ and $r_2,r_3\geq r-d+1$ 
will have $f_d(P_{\cal{F}}) =f_{2;d}(r)$ in this case. (ii) Otherwise,
if $0\leq d < 2(r-1)/3$, then the Minkowski sum  $P_{\cal{F}} = \Delta_{F_1} + \Delta_{F_2}$
with $f_d(P_{\cal{F}}) =f_{2;d}(r)$ is more ``unique'' in the following sense:

If $r = 3q$, then $F_1$ and $F_2$ must satisfy
\[
(r_1,r_2,r_3) \in \{ (q,q,q), (q+1,q,q-1), (q+1,q-1,q)\}.
\]

If $r = 3q+1$, then $F_1$ and $F_2$ must satisfy 
$(r_1,r_2,r_3) = (q+1,q,q)$.

If $r = 3q+2$, then $F_1$ and $F_2$ must satisfy
\[
(r_1,r_2,r_3) \in \{ (q+2,q,q), (q+1,q+1,q), (q+1,q,q+1)\}.
\]

Directly from Theorem~\ref{thm:max-2d-faces} we have the following.
\begin{corollary}
\label{cor:assymp-d}
For a fixed $d$ we have
\[
\lim_{r\rightarrow\infty}\frac{f_{2;d}(r)}{r^{d+2}} = \frac{2(3^{d+1}-2^{d+1})}{3^{d+1}(d+2)!}.
\]
\end{corollary}

\subsection{Second Example}

Consider the Minkowski sum $P_{\cal{F}} = \Delta_{F_1} + \Delta_{F_2}$ of two sets ${\cal{F}} = (F_1,F_2)$,
where $r = r_1+r_2+r_3$, $r_1 = |F_1\cap F_2|$, $r_2 = |F_1\setminus F_2|$ and $r_3 = |F_2\setminus F_1|$
as before. Directly by Corollary~\ref{cor:k=2} we have
that $f_0 := f_0(P_{\cal{F}})$, the number of vertices of $P_{\cal{F}}$,
is given by $f_0 = r_1r_2 + r_1r_3 + r_2r_3 + r_1$, something which we saw
alternatively in the last example of Section~\ref{sec:flag-poly}. If $F_1$ and $F_2$ are disjoint,
then $P_{\cal{F}} = \Delta_{F_1}\times\Delta_{F_2}$ (see Proposition 1.2 in~\cite{Minkowski} for
a more general result) which is a simple polytope by the simplicity of
any simplex. Hence we will assume that $F_1$ and $F_2$ have a nonempty intersection
so $r_1\geq 1$, in which case the dimension of $P_{\cal{F}}$ is precisely $r-1$ 
(see~\cite{Minkowski} for more detailed discussion).

For a polynomial $f(\tilde{x})$, denote by $[{\tilde{x}}^{\tilde{a}}]f(\tilde{x})$
the coefficient for ${\tilde{x}}^{\tilde{a}}$ in $f(\tilde{x})$.
Note that the number of $(0,d)$-chains in $P_{\cal{F}} = \Delta_{F_1} + \Delta_{F_2}$
is then given by $[x_2^d]{\tilde{f}}^2_{P_{\cal{F}}}(x_1,x_2)$.
Note that when $F_1$ and $F_2$ are not disjoint, then each vertex of 
$P_{\cal{F}}$ is contained in a minimum of $\binom{r-1}{d}$
$d$-faces since the degree of each vertex is at least $r-1$ and $d$ of these yield
a face. 
\begin{observation}
\label{obs:(0,r-2)simple}
For $P_{\cal{F}} = \Delta_{F_1} + \Delta_{F_2}$ where $F_1$ and $F_2$ have a
nonempty intersection, $P_{\cal{F}}$ is simple iff 
\[
[x_2^d]{\tilde{f}}^2_{P_{\cal{F}}}(x_1,x_2) = \binom{r-1}{d}f_0 = \binom{r-1}{d}(r_1r_2 + r_1r_3 + r_2r_3 + r_1).
\]
\end{observation}
Since $[x_2^d]{\tilde{f}}^2_{P_{\cal{F}}}(x_1,x_2) = [x_2^d]{\tilde{f}}^2_{P_{\cal{F}}}(0,x_2)$,
and in addition ${\tilde{f}}^2_{\Delta_{r_i-1}}(0,x) = r_i(x+1)^{r_i-1}$ and
${\tilde{f}}^1_{\Delta_{r_i-1}}(x) = ((x+1)^{r_i} - 1)/x$,
we obtain by Corollary~\ref{cor:2-chain} that 
\begin{eqnarray*}
f_{(0,d)}(P_{\cal{F}}) & = & [x^d]{\tilde{f}}^2_{P_{\cal{F}}}(0,x) \\
                       & = & [x^d]\left( (r_1r_2 + r_1r_3 + r_2r_3 + r_1)(x+1)^{r-1}
                                               + r_2r_3(x+1)^{r-2} - r_2r_3(x+1)^{r_2+r_3-1}\right),
\end{eqnarray*}
and hence we have the following.
\begin{corollary}
\label{cor:(0,d)chains}
The number $f_{(0,d)}(P_{\cal{F}})$ of $(0,d)$-chains in $P_{\cal{F}} = \Delta_{F_1} + \Delta_{F_2}$
is given by 
\[
 f_{(0,d)}(P_{\cal{F}}) = (r_1r_2 + r_1r_3 + r_2r_3 + r_1)\binom{r-1}{d}
                       + r_2r_3\binom{r-2}{d} - r_2r_3\binom{r_2+r_3-1}{d}.
\]
\end{corollary}
By Observation~\ref{obs:(0,r-2)simple} and Corollary~\ref{cor:(0,d)chains}
it is immediate that the Minkowski sum $P_{\cal{F}} = \Delta_{F_1} + \Delta_{F_2}$
is simple iff $r_2=0$, $r_3=0$ or $r_1 = 1$. 

The other 
question {\em ``when is $P_{\cal{F}} = \Delta_{F_1} + \Delta_{F_2}$ 
farthest from being simple?''} is more interesting. To make this question precise
and address it properly, we define the {\em discrepancy} $\delta_{(0,d)}(P_{\cal{F}})$ 
as the difference between $f_{(0,d)}(P_{\cal{F}})$ and the minimum number $\binom{r-1}{d}f_0$
of $(0,d)$-chains in simple polytope of the same dimension and on the same
number $f_0 = r_1r_2 + r_1r_3 + r_2r_3 + r_1$ of vertices,
\[
\delta_{(0,d)}(P_{\cal{F}}) := f_{(0,d)}(P_{\cal{F}}) 
- \binom{r-1}{d}f_0 = r_2r_3\left(\binom{r-2}{d} - \binom{r_2+r_3-1}{d}\right).
\]
The discrepancy $\delta_{(0,d)}(P_{\cal{F}})$ can be thought of as a 
measure (in terms of $d$) on how far from simple the polytope 
$P_{\cal{F}}$ is. Hence, for a fixed $r\in\nats$ we like to determine
what kind of Minkowski sum $P_{\cal{F}}$ of two simplices
will yield the largest discrepancy. 
This is only of interest if $r_2r_3\neq 0$ and $r_1\geq 1$,
\begin{equation}
\label{eqn:discrep}
\delta_{2;(0,d)}(r) := \max\{\delta_{(0,d)}(P_{\cal{F}}) : P_{\cal{F}} = \Delta_{F_1} + \Delta_{F_2}, \ \ 
r_1 + r_2 + r_3 = r, \ \ r_1\geq 1\}.
\end{equation}
Since $r_1\geq 1$ then $r_2+r_3\leq r-1$, so to determine $\delta_{2;(0,d)}(r)$ it suffices
to compute $\max(f)$ where $f(x,y)$ is given by 
\[
f(x,y) = \left\{ \begin{array}{ll}
                   xy\binom{r-2}{d} & \mbox{ if } 0\leq x+y\leq d, \ \ x,y\geq 0, \\ 
                   xy\left(\binom{r-2}{d} - \binom{x+y-1}{d}\right) & \mbox{ if } d+1\leq x+y\leq r-1, \ \ x,y\geq 0.
                 \end{array}         
         \right. 
\]
Considering $f(x,y)$ in integers $x$ and $y$ on each line $x+y = c$ 
where $c$ is a constant, then clearly 
the maximum of $f(x,y)$ is reached on each of these lines when
$|x-y|\in\{0,1\}$. In particular, when $0\leq x+y\leq d$ the maximum of $f(x,y)$ is 
$\lfloor d/2\rfloor\lceil d/2\rceil\binom{r-2}{d}$. For the other range 
$d+1\leq x+y\leq r-1$, the more interesting one (since we like to consider
$d$ fixed and $r$ much larger than $d$),
we can by symmetry assume that $x\leq y$ and write $y = x +\epsilon$
where $\epsilon \in\{0,1\}$. 
Then by expressing $f(x,y) = f(x,x+\epsilon)$ in terms of $z = x+y-1 = 2x+\epsilon -1$,
we have that $\delta_{2;(0,d)}(r) = \max(g)$ where
$g(z) = \frac{1}{4}(z+1-\epsilon)(z+1+\epsilon)\left(\binom{r-2}{d} - \binom{z}{d}\right)$
and $d\leq z\leq r-2$. For the moment we consider
real $z$ in the mentioned interval, so by the intermediate value theorem
for real continuous functions, $\max(g)$ is obtained at $g(z^*)$ where
$z-1\leq z^*\leq z$ and the difference $\Delta g(z-1) = g(z) - g(z-1) = 0$, or
\begin{equation}
\label{eqn:Dg=0}
(z+1-\epsilon)(z+1+\epsilon)\left(\binom{r-2}{d} - \binom{z}{d}\right)=
(z-\epsilon)(z+\epsilon)\left(\binom{r-2}{d} - \binom{z-1}{d}\right).
\end{equation}
We will first estimate asymptotically the solution
in $z$ to (\ref{eqn:Dg=0}), and hence also $z^*\in [z-1,z]$, in terms of 
$d$ and $r$ and interpret the corresponding $\delta_{2;(0,d)}(r) = \max(g)$.

Note that (\ref{eqn:Dg=0}) can be rewritten as 
\begin{equation}
\label{eqn:prodh}
\prod_{i=1}^d\frac{r-i-1}{z-i+1} = \frac{(d+2)z + 1 + \epsilon^2d/z}{2z+1}.
\end{equation}
Considering solely the asymptotic values as $r$ tends to infinity, the solution $z$
of (\ref{eqn:prodh}) must also tend to infinity. Hence, as the right hand side (RHS)
of (\ref{eqn:prodh}) tends to $(d+2)/2$ as $r$ (and hence $z$ as well) goes
to infinity, and the left hand side (LHS) of (\ref{eqn:prodh})
is asymptotically $(r/z)^d$ as $r$ tends to infinity, we have 
for large $r$ that $z\sim \sqrt[d]{2/(d+2)}r$. Hence, $\max(g)$ is obtained at $g(z^*)$ 
where $z^*\sim \sqrt[d]{2/(d+2)}r$ as $r$ tends to infinity. From this we deduce
the following asymptotic result.
\begin{proposition}
\label{prp:discrep-asymp}
If a Minkowski sum $P_{\cal{F}} = \Delta_{F_1} + \Delta_{F_2}$ 
where $F_1, F_2\subseteq \{1,\ldots,r\}$ and $r_1\geq 1$ has the 
maximum discrepancy $\delta_{2;(0,d)}(P_{\cal{F}}) = \delta_{2;(0,d)}(r)$ between 
the number of $(0,d)$-chains
of faces of $P_{\cal{F}}$ and the minimum number
$\binom{r-1}{d}f_0$ of $(0,d)$-chains in a simple 
polytope of the same dimension and with the same
number $f_0$ of vertices,
then for all fixed $d$ and $i\in \{2,3\}$ we have
\begin{eqnarray*}
r_1 & \sim & \left(1 - \sqrt[d]{\frac{2}{d+2}}\right)r, \\
r_i & \sim & \frac{1}{2}\sqrt[d]{\frac{2}{d+2}}r, 
\end{eqnarray*}
as $r$ tends to infinity, and hence
\[
\lim_{r\rightarrow\infty}\frac{\delta_{2;(0,d)}(r)}{r^{d+2}} 
= \frac{2^{2/d - 2}}{(d+2)^{2/d-1}(d-1)!}.
\]
\end{proposition}

We will now determine exactly which $r_1$, $r_2$ and $r_3$ will yield the maximum
discrepancy. In order to do that we 
determine concrete and tight upper and lower bounds for the solution $z$ of 
(\ref{eqn:prodh}) to pinpoint the corresponding integer $z^*$ at which
the maximum discrepancy $\delta_{2;(0,d)}(r) = \max(g) = g(z^*)$ for the 
Minkowski sum $P_{\cal{F}} = \Delta_{F_1} + \Delta_{F_2}$ is attained.
The remainder of this section will be spent computing this mentioned
$z^*$ and interpreting the corresponding discrepancy for the Minkowski sum.

First we centralize the LHS of (\ref{eqn:prodh}) by letting $R = r - (d+3)/2$
and $w = z - (d-1)/2$. In this way the LHS of (\ref{eqn:prodh}) becomes 
\begin{equation}
\label{eqn:LHS}
\prod_{i=1}^{d/2}\frac{R^2 - (i-1/2)^2}{w^2 - (i-1/2)^2}\ \ \mbox{ for even $d$, and }\ \ 
\frac{R}{w}\prod_{i=1}^{(d-1)/2}\frac{R^2 - i^2}{y^2 - i^2}\ \ \mbox{ for odd $d$.}
\end{equation}
Since $w<R$ are fixed then $i\mapsto h(i) := (R^2-i^2)/(w^2-i^2)$ is increasing and the 
RHS of (\ref{eqn:prodh}) is at most
$(d+2)/2$, we have that $(R/w)^d < (d+2)/2$ and therefore we obtain a concrete lower bound
$L(d;r)$ for $w$, and hence also for $z$, given by
\begin{equation}
\label{eqn:zlower}
z > L(d;r) := \sqrt[d]{\frac{2}{d+2}}(r-2) + \left(1-\sqrt[d]{\frac{2}{d+2}}\right)\frac{d-1}{2}.
\end{equation}
Incidentally, we will use this lower bound $L(d;r)$ from (\ref{eqn:zlower}) to obtain an
upper bound for $z$. By (\ref{eqn:prodh}), its LHS expressed in (\ref{eqn:LHS}) is at least 
$q(z) = ((d+2)z+1)/(2z+1)$. Since $z\mapsto q(z)$ is increasing, the RHS of (\ref{eqn:prodh}) and hence 
also (\ref{eqn:LHS}) is greater than $q(L(d;r))$. 
Regardless of the parity of $d$, the expression in (\ref{eqn:LHS}) is less than $h((d-1)/2)^{d/2}$ 
and hence we obtain an upper bound for $w$, and hence also $z$, determined by
$q(L(d;r))<h((d-1)/2)^{d/2}$ or 
\begin{equation}
\label{eqn:zupper}
w < cR\sqrt{1 + \left(\frac{d-1}{2}\right)^2\left(\frac{1}{c^2}-1\right)\frac{1}{R^2}},
\end{equation}
where $c = 1/\sqrt[d]{q(L(d;r))}$. 

By the mean value theorem for the function $t\mapsto \sqrt[d]{t}$, there is an 
$\alpha \in \left[\frac{2}{d+2},\frac{1}{q(L(d;r))}\right] = \left[\frac{2}{d+2}, c^d \right]$ such that 
\[
0 < c - \sqrt[d]{\frac{2}{d+2}} = 
\frac{1}{{\alpha}^{1-1/d}(d+2)((d+2)L(d;r)+1)} < \frac{1}{2(d+2)(r-2)},
\]
and since $\sqrt{1+t} < 1 + t/2$ for positive $t$, (\ref{eqn:zupper}) implies 
\begin{equation}
\label{eqn:zzupper}
z < cR + \frac{d-1}{2} + 
  \frac{1}{2}\left(\frac{d-1}{2}\right)^2\left(\frac{1}{c}-c\right)\frac{1}{R}.
\end{equation}
By the bounds of $c$ and definition of $L(d;r)$ we have further that
\begin{eqnarray*}
cR + \frac{d-1}{2} & < & \left(\sqrt[d]{\frac{2}{d+2}}  
                   + \frac{1}{2(d+2)(r-2)}\right)R + \frac{d-1}{2} \\
                   & = & L(d;r) + \frac{R}{2(d+2)(r-2)} \\
                   & < & L(d;r) + \frac{1}{2(d+2)},
\end{eqnarray*}
and hence (\ref{eqn:zzupper}) implies 
\begin{equation}
\label{eqn:z3upper}
z < L(d;r)  + \frac{1}{2(d+2)} + 
  \frac{1}{2}\left(\frac{d-1}{2}\right)^2\left(\frac{1}{c}-c\right)\frac{1}{R}. 
\end{equation}
Since $t\mapsto 1/t-t$ is decreasing we have further
\[
\frac{1}{c}-c < \sqrt[d]{\frac{d+2}{2}} + \sqrt[d]{\frac{2}{d+2}} < \frac{2\log(d+2)}{d},
\]
and so (\ref{eqn:z3upper}) implies the following.
\begin{proposition}
\label{prp:tight}
If $L(d;r)$ is as defined in (\ref{eqn:zlower}), then the real number solution $z$ to 
(\ref{eqn:Dg=0}) satisfies 
\[
0 < z-L(d;r) < \frac{1}{2(d+2)} + \frac{d\log(d+2)}{2(2r-d-3)}
\]
for all natural numbers
$d$ and $r$ with $d < r$.
\end{proposition}

\noindent{\sc Note:} (i) By Proposition~\ref{prp:tight} we have  
for any $\varepsilon > 0$ there is an integer $d'_{\varepsilon}$ and 
a function $t\mapsto r'_{\varepsilon}(t)$
such that for any $d\geq d'_{\varepsilon}$ and $r\geq r'_{\varepsilon}(d)$ then
$0 < z-L(d;r) < \varepsilon$. (ii) Also, as 
$(1-\sqrt[d]{2/(d+2)})(d-1)/2 = \log(d/2)/2 + O(\log^2d/d)$
we immediately obtain that the real number solution $z$ to (\ref{eqn:Dg=0}) satisfies
\[
z = \sqrt[d]{\frac{2}{d+2}}(r-2) + \frac{\log(d/2)}{2} 
+ O\left(\frac{\log^2(d)}{d}\right)
+ O\left(\frac{d\log(d)}{r}\right),
\]
for large $d$ and $r$ with $d < r$.

\vspace{3 mm}

By the above Proposition~\ref{prp:tight} we now have the following result.
\begin{theorem}
\label{thm:discrep}
If a Minkowski sum $P_{\cal{F}} = \Delta_{F_1} + \Delta_{F_2}$ 
where $F_1, F_2\subseteq \{1,\ldots,r\}$ and $r_1\geq 1$ has the 
maximum discrepancy $\delta_{2;(0,d)}(P_{\cal{F}}) = \delta_{2;(0,d)}(r)$ between 
the number of $(0,d)$-chains
of faces of $P_{\cal{F}}$ and the minimum number
$\binom{r-1}{d}f_0$ of $(0,d)$-chains in a simple 
polytope of the same dimension and with the same
number $f_0$ of vertices,
then we have $|r_2-r_3|\leq 1$ and
\[
\left\lceil  L(d;r)\right\rceil \leq r_2+r_3 \leq
\left\lfloor L(d;r) + \frac{2d+5}{2d+4} + \frac{d\log(d+2)}{2(2r-d-3)}\right\rfloor,
\]
where $L(d;r)$ is as defined in (\ref{eqn:zlower}).
\end{theorem}

Since the ring $\ints$ of the integers is a UFD, the proof of
the following is an easy exercise.
\begin{claim}
\label{clm:irrational}
For any integer $d\geq 2$ the number $\sqrt[d]{\frac{2}{d+2}}$ is irrational.
\end{claim}
By the above Claim~\ref{clm:irrational} and the Equidistribution Theorem
by Hermann Weyl~\cite[p.~107]{Stein-Shak}, we have the following.
\begin{observation}
\label{obs:dense}
For every fixed integer $d\geq 2$ the set $\{L(d;r)\bmod 1 : r\in\nats\}$ is
uniformly distributed and hence a dense subset of the unit interval $[0,1]$ 
with respect to the usual topology of the reals.
\end{observation}
Consider the upper bound of $r_2+r_3$ in Theorem~\ref{thm:discrep}.
For every fixed $d\geq 1$ we see that if $r > d\log(d+2)$ then
the real interval in which the integer $r_2+r_3$ is contained in
is of length strictly less than 2 and hence there are at most two possibilities
for the integer $r_2+r_3$. Observation~\ref{obs:dense} now implies
there are infinitely many $r > d\log(d+2)$ for which the integer
$r_2+r_3$ is uniquely determined by Theorem~\ref{thm:discrep}.
In this case we know exactly when the maximum discrepancy occurs.
\begin{corollary}
\label{cor:exact}
With the setup from Theorem~\ref{thm:discrep}, if 
$\delta_{2;(0,d)}(P_{\cal{F}}) = \delta_{2;(0,d)}(r)$, then for all fixed $d\geq 1$,
$r > d\log(d+2)$ and $i\in \{2,3\}$, we have
\begin{eqnarray*}
|r_2 - r_3| & \leq & 1, \\
|r_2 + r_3 - \left\lceil L(d;r)\right\rceil| & \leq & 1.
\end{eqnarray*}
By symmetry we may assume $r_2\leq r_3$ to hold, in which case
we have further there are infinitely many $r$ such that 
the maximum discrepancy $\delta_{2;(0,d)}(r)$ occurs exactly at
\begin{eqnarray*}
r_1 & = & r - \left\lfloor \left\lceil L(d;r)\right\rceil/2\right\rfloor
            - \left\lceil \left\lceil L(d;r)\right\rceil/2\right\rceil, \\
r_2 & = & \left\lfloor \left\lceil L(d;r)\right\rceil/2\right\rfloor, \\
r_3 & = & \left\lceil \left\lceil L(d;r)\right\rceil/2\right\rceil.
\end{eqnarray*}
\end{corollary}
\noindent{\sc Remark:}
Although $r_2+r_3$ is for all $r>d\log(d+2)$ uniquely determined up to
parity, we have by the uniform distribution property of $L(d;r)\bmod 1$ in 
Observation~\ref{obs:dense} that for every fixed $d$, {\em ``most''} $r>d\log(d+2)$
will determined $r_1$, $r_2$ and $r_3$ uniquely as stated in the latter
display in Corollary~\ref{cor:exact}. That is to say, the probability
of picking $r$ that uniquely determines $r_1$, $r_2$ and $r_3$ can be
made arbitrarily high by choosing large enough $d$ and large enough
$r>d\log(d+2)$.

\subsection*{Acknowledgments}

The author would like to thank James F.~Lawrence
for helpful discussions regarding the theory of polytopes
in general. Also, thanks to Bernd Sturmfels for
introducing this investigation of Minkowski sum of standard simplices.
Last but not least, sincere thanks to Walter Morris for his
encouragement, interest and helpful discussions. 
Figure~\ref{fig:contracted-cube}
is the same that appears in~\cite{Minkowski} and was drawn by Walter Morris.

\flushright{\today}
\end{document}